\documentclass[letterpaper, 11pt]{article} 
\usepackage[top=1in, bottom=1in, left=1in, right=1in]{geometry}

\usepackage{amsmath,amssymb,amsthm}
\newtheorem{Thm}{Theorem}

\newtheorem{Def}{Definition}
\usepackage{epsfig}
\usepackage{arydshln}
\usepackage{subcaption}

\usepackage{comment}

\newcommand{\bfm}[1]{\mbox{\boldmath ${#1}$}}
\newcommand{\nonum}{\nonumber \\}
\newcommand\eq[1] {(\ref{#1})} 
\newcommand{\beqa}{\begin{eqnarray}}
\newcommand{\eeqa}[1]{\label{#1}\end{eqnarray}}
\newcommand{\beq}{\begin{equation}}
\newcommand{\eeq}[1]{\label{#1}\end{equation}}
\newcommand{\Grad}{\nabla}

\newcommand{\Md}{\partial}

\newcommand{\Ga}{\alpha}

\newcommand{\Ge}{\epsilon}

\newcommand{\Gf}{\phi}

\newcommand{\Gl}{\lambda}

\newcommand{\GG}{\Gamma}

\newcommand{\GO}{\Omega}

\newcommand{\BGs}{\bfm\sigma}

\newcommand{\CO}{{\cal O}}

\newcommand{\bpm}{\begin{pmatrix}}
\newcommand{\epm}{\end{pmatrix}}
\newcommand\fig[1] {{\rm Figure}~\ref{fig:#1}}

\newcommand\labfig[1] {\label{fig:#1}}
\newcommand\sect[1] {\ref{sect:#1}}
\newcommand\labsect[1] {\label{sect:#1}}
\usepackage{url}
\usepackage{hyperref}

\def\Bf{{\bf f}}
\def\Bg{{\bf g}}

\def\Bt{{\bf t}}
\def\Bu{{\bf u}}
\def\Bv{{\bf v}}
\def\Bw{{\bf w}}
\def\Bx{{\bf x}}
\def\By{{\bf y}}

\def\BR{{\bf R}}

\usepackage{graphics}
\usepackage{amssymb}

\usepackage{color}

\newcommand{\cb}{\color{black}}
\newcommand{\cred}{\color{black}}
\newcommand{\crev}{\color{black}}
\newcommand{\cn}{\color{black}}
\newcommand{\cg}{\color{black}}
\newcommand{\crevv}{\color{black}}

\begin{document}

\title{The obstacle problem in masonry structures and cable nets}

\author{Ada Amendola$^{1}$,
  Ornella Mattei$^{2}$,
  Graeme W. Milton$^{3}$,
Pierre Seppecher$^{4}$\\
\medskip\\
\small$^{1}$Department of Civil Engineering, University of Salerno, Fisciano (SA), Italy. \\
  \small$^{2}$Department of Mathematics, San Francisco State University, CA, USA. \\
 \small $^{3}$Department of Mathematics, University of Utah, Salt Lake City, UT, USA. \\
\small$^{4}$Institut de Math\'ematiques de Toulon,  Universit\'e de Toulon et du Var, France.}



\date{}
\maketitle
\begin{abstract}
 \noindent
 \crev We consider the problem of  finding a \cn \cb net \cn \crev that supports prescribed forces applied at prescribed points, yet avoids certain
obstacles, with all the elements of the \cn \cb net \cn \crev under compression (strut net) or under tension (cable web). In the case of masonry structures, for instance, this consists in finding a strut net that supports the forces, \cn \cb is \cn \crevv contained \cn within the physical structure, and avoids regions that may be not accessible due, for instance, to the presence of holes. We solve such a problem in the two-dimensional case, where the prescribed forces are applied at the vertices of a convex polygon, and we treat the cases 
of both single and multiple obstacles. \cn
  By approximating the obstacles by polygonal regions, the task reduces to identifying the feasible \cb domain \cn  in a linear programming problem.
  For a single obstacle we show how the
  region $\GG$ available to the obstacle can be enlarged as much as possible in the sense that there is no other strut net, having a region
  $\GG'$ available to the obstacle with $\GG\subset\GG'$. The case where some of the forces are reactive, unprescribed but reacting to the
  other prescribed forces, is also treated. It again reduces to identifying the feasible \cb domain \cn  in a linear programming problem. Finally,
  one may allow a subset of the reactive forces to each act not at a prescribed point, but rather at any point on a prescribed line
  segment. Then the task reduces to identifying the feasible \cb domain \cn  in a quadratic programming problem. 

\end{abstract}





\maketitle

\section{Introduction} \labsect{1}
The master `safe' theorem of a masonry arch, formulated in a famous work by Jacques Heyman \cite{heyman1995} states that \crevv the collapse of such a structure \cn\cg does not \cn\crevv occur 
if a line of thrust of the given external loads can be fit within the boundaries of the arch
(see also, \cite{Alexakis:2014:LEA,Alexakis:2015:LEA,como2017}).
\cn
The line of thrust coincides with the Hooke's inverted chain of the given forces, as shown, e.g., by the Italian engineer Poleni in a study dated back to 1784, which was aimed at assessing the stability of St. Peter's dome in \cg the \cn Vatican \cite{poleni}. Such a construction has been the foundation of the beautiful work of the Catalan architect Antoni Gaud\'{i} on funicular shapes of masonry structures \cite{gaudi,huerta}.
The use of truss networks within discrete element approaches to no-compression or no-tension bodies is frequent in the literature
(refer to \cite{Bouchitte:2019:FCW,Milton:2017:SFI,odwyer1999,block2007,marmo2017}
and references therein).
Truss models are convenient, e.g., to tackle the equilibrium problem of masonry structures, described as no-tension bodies \cite{Giaquinta:1988:RMS,delpiero1998}, and to assess the compatibility of external loads through graphical constructions,  numerical methods, and hanging models
\cite{odwyer1999,block2007,marmo2017,lucchesi2008,fortunato2010,defaveri2013}.
Polyhedral Airy stress functions have been employed in two-dimensional
problems, since such functions permit one to conveniently describe force networks through scalar potentials  \cite{fraternali2010,Fraternali:2014:CBF}. 
\crevv
The no-tension model of masonry by Heyman does not 
account for buckling analysis, since the material is supposed to behave rigidly up to collapse, with infinitely large compression strength \cite{heyman1995}. 
Upon accounting for a geometrically nonlinear elastic response and a finite strength in compression of the masonry, it has been shown that second-order effects on the collapse load of circular arches tend to be negligible \cn\cg when the slenderness ratio of the arch radius over its thickness is sufficiently small \cite{Zani2009}.
\cn

A recent study has investigated the existence problem of systems of
tensile forces that support given sets of nodal forces, coming to the
conclusion that if such `cable webs' exist then the given forces are
supported by the complete web connecting pairwise their points of
application (see \cite{Bouchitte:2019:FCW}, Theorem 1.1).  In
two-dimensions, and with forces at the vertices of a convex polygon, the use of polyhedral Airy stress functions shows
alternatively that if such `cable webs' exist then the given forces are  supported by an open web with no internal loops (see
\cite{Milton:2017:SFI}, Theorem 1).  
The present work generalizes the above problem to the case  of compression-only force networks (or `strut nets')  in two-dimensions, \cn which support given sets of nodal forces, and are able to avoid `obstacles' representing regions not accessible to the force network (e.g., holes or inclusions). 
We begin by deriving a result showing that a strut net supporting $q$ forces located strictly inside the convex hull of
the points where forces are applied and  avoiding $p$
obstacles can be replaced by a strut net
with at most $p+q$ elementary loops. (Section \sect{2}). From then
onwards we restrict our attention to sets of forces at the vertices of
a convex polygon ($q=0$). In Section \sect{3} we provide a condition for the non-existence of strut nets that avoid a given obstacle, when such an obstacle intersects the open strut net associated to the given forces,
and extends outside the convex hull of their points of application. We also provide an algorithm to establish if there
exists a strut net avoiding multiple obstacles, and if so to
construct one. Additionally, for a single obstacle,  we show how the region available to the
obstacle can be enlarged as much as possible  (in a sense to be made more precise). 
Section \sect{4} treats the case of multiple obstacles with the
inclusion of reactive forces.  Numerical applications are then presented in
Section \sect{5}, by first dealing with a test problem characterized by
a single active force, and then by treating a series of examples that
analyze strut nets associated with the statics of masonry arches. We
end in Section \sect{6} by drawing concluding remarks and directions for future work.

\section{Reducing the complexity of a strut net avoiding a given set of obstacles} \labsect{2}
\setcounter{equation}{0}

\crev To begin, we consider a set of points $\Bx_1,\Bx_2,\ldots,\Bx_n$ at the vertices of \cn \cb a convex polygon $\GO$, \cn \cb with external \cn \crev
forces $\Bt_1,\Bt_2,\ldots,\Bt_n$ acting \cn \cb on \cn \crev them. We assume that the points are numbered anti clockwise and we adopt the following convention : for any $i>n$, $\Bx_i:=\Bx_{i-n}$ and 
$\Bx_0:=\Bx_{n}$. \cn
Letting \cb
$ \BR_\perp=\bpm 0 & -1 \\ 1 & 0 \epm $ \cn
denote the matrix for a $90^\circ$ rotation, we assume that, for all $j$ with
$1\leq j \leq n$ and for all $i$ with $j \leq i \leq j+n-1$,
\beq  \sum_{k=j}^{i}(\Bx_k-\Bx_j)\cdot[\BR_\perp\Bt_k]\geq 0.
\eeq{0.1b}
\cn This assumption ensures that there exists a truss structure supporting
these forces with all the struts in the truss structure under compression \cite{Milton:2017:SFI}.
The constraint \eq{0.1b} has a physical interpretation: the net anticlockwise torque around
the point $\Bx_j$ of the forces $\Bt_k$ summed over any number of consecutive points clockwise past the point
$\Bx_j$ is non-positive. If \eq{0.1b} is satisfied, the forces are supported by an open strut net with no internal loops, \crev in other words, there are no polygons whose edges are the
\cn \cb struts of the net. \cn

A more general condition, applicable to any set of forces, either in two-dimensions or three-dimensions, is that if a strut net exists
supporting the forces, then they will also be supported by a strut net \crev connecting pairwise the points at which the forces are applied \cn\cite{Bouchitte:2019:FCW}. This leads to a linear
programming problem for determining if a set of forces can be supported. 

\cb Condition \eq{0.1b} is derived using Airy stress functions. \cn  Indeed, for two-dimensional elasticity it is well known that, in the absence of body forces in a simply connected
region $\GO$, the \crev divergence-free  stress field $\BGs(\Bx)$ \cn
can be represented in terms of the Airy stress function $\phi(\Bx)$:
\beq \BGs(\Bx)=
\bpm \frac{\Md^2\phi(\Bx)}{\Md^2 x_2^2} &-\frac{\Md^2\phi(\Bx)}{\Md x_1 \Md x_2} \\ -\frac{\Md^2\phi(\Bx)}{\Md x_1 \Md x_2}& \frac{\Md^2\phi(\Bx)}{\Md^2 x_1^2}\epm =
\BR_\perp^T\Grad\Grad\phi(\Bx)\BR_\perp, \eeq{0.1}
in which $\BR_\perp^T=-\BR_\perp$ is the transpose of $\BR_\perp$. 
Since $\BGs(\Bx)$ is negative semidefinite for all $\Bx$, we
see that $\Grad\Grad\phi(\Bx)$ is negative semidefinite for all $\Bx$. 
This implies that $\phi(\Bx)$ is a concave function in any simply-connected two-dimensional region under compression, that may
have subregions with zero stress \cite{Giaquinta:1988:RMS}.
Note that, when $\BGs(\Bx)$ is zero in a region, as it is between the struts in a net, then $\Grad\Grad\phi(\Bx)=0$
which by \eq{0.1} implies the Airy stress function is a linear function of $\Bx$  in this region.
Thus, any strut net under compression that supports
forces at the vertices of a convex polygon will have an associated Airy stress function \cb whose graph \cn is a concave polyhedron.
The magnitudes of the slope discontinuities can be connected to the compression in the associated struts \cite{Fraternali:2014:CBF, Fraternali:2002:LSM}.

\crev If the strut net associated with a concave polyhedron has internal loops, we can replace  the latter
by a simpler concave polyhedron (supporting the same forces),  by using the tangent planes \cn \cb to the Airy stress function \cn \crev at the boundary of the polygon $\GO$
\cite{Milton:2017:SFI}. \cn
The strut net associated with
this latter concave polyhedron will be an open strut net with no internal loops. \crev To prove this, let us start by noticing that, since the strut net \cn supports under compression the forces $\Bt_i$ at the points $\Bx_i$, it can support, for some $\epsilon>0$, the same    
forces at points 
\beq \tilde{\Bx}_{i}=\Bx_i-\Ge\Bt_i, \eeq{0.4}
the convex hull of which we call $\tilde \Omega$. To do this, one extends the strut net
by attaching $n$ short struts of length $\Ge>0$ between $\Bx_i$ and $\tilde{\Bx}_{i}$, $i=1,2,\ldots,n$.
The first step is to determine the Airy stress function $\phi(\Bx)$ in the polygonal ring bounded on one side
by the polygon joining the points $\Bx_i$, and on the other side by the polygon joining the points $\tilde{\Bx}_{i}$.
When $\Ge$ is sufficiently small there are no struts inside the quadrilateral \crev $\GO_i$ \cn
with vertices $\tilde{\Bx}_{i-1}$, $\Bx_{i-1}$, $\Bx_i$, $\tilde{\Bx}_{i}$
and since the stress vanishes there, the Airy stress function $\phi(\Bx)$ inside that quadrilateral must be a linear function
$L_i(\Bx)$. Due to the stress in the wires (see, for example, \cite{Milton:2017:SFI}), these satisfy
\beq \BR_\perp(\Grad L_{i+1}-\Grad L_{i})=-\Bt_i,\quad i=1,2,\ldots, n, \eeq{A.1}
and, for all $i,j$,
\beq L_i\geq L_j\quad\text{on }\GO_j. \eeq{A.2}
Equation \eq{A.1} and the constraints
\beq L_i(\Bx_i)=L_{i+1}(\Bx_i) \eeq{A.2a}
determine the sequence $(L_i)_{i=1}^n$ up to the addition of a global linear function.
Uniqueness can be ensured by assuming, for instance, that $L_1(\Bx)=0$. The
condition \eq{A.2} is equivalent to \eq{0.1b} \cite{Milton:2017:SFI}.

\cb The \cn \crev stress measure $\BGs(\Bx)$ associated to the Airy function $\phi(\Bx)$, that is defined in the sense of distributions by \eq{0.1}, is non-positive, concentrated along a finite union of segments, divergence-free on $\Omega$ and, owing to the construction of $(L_i)_{i=1}^n$, \cn\cb the divergence in the sense of distributions, $div(\sigma_{|\Omega})$,
coincides \cn with the discrete measure $\sum_{i=1}^n \Bt_i \delta_{\Bx_i}$ on $\partial \Omega$. \cb Let us introduce the following definition: \cn
\begin{Def} A ``strut net function" is a concave piecewise linear function $\phi(\Bx)$ on $\widetilde \Omega$ coinciding with $L_i(\Bx)$ on each $\Omega_i$. 
\end{Def}
\noindent \cb We \cn see that, in particular, the envelope of the tangent planes,
\beq \phi_0(\Bx)\equiv\min_i L_i(\Bx), \eeq{0.9}
is a strut net function: we call it the ``open strut net function" as it is associated with a strut net having no loops. \cn Note also that concavity implies that any strut net function $\phi(\Bx)$ satisfies, for all $i$,  $\phi(\Bx)\leq L_i(\Bx)$, and hence
that  $\phi(\Bx)\leq \phi_0(\Bx)$.
\subsection{Forces at points inside the convex hull of all application points and multiple obstacles}

\crev Consider a strut net containing  a convex elementary polygonal loop \cn \cb within which there are no obstacles. By \cn \crev
elementary we mean a loop with no interior struts. Then, the net forces (due to adjoining struts or imposed
forces) act on the vertices of the loop as in \fig{0}(a), and the polygonal loop can be replaced by an open strut net, \crev as in \fig{0}(b).
In this way, if there are many internal loops in the strut net, we decrease the
number of elementary loops by one. This reduction can be continued \cn until any remaining elementary loop encloses one or more obstacles
or has an imposed
force acting outwards from a vertex of the polygonal loop, so that the
loop is not convex.  

Then, using the fact that any obstacle intersecting
the boundary of the convex hull of all points where forces are applied
can never have a loop of struts under compression surrounding it,
we establish:

\begin{Thm}
  A strut net that has $q$ imposed forces at points that lie strictly inside the convex hull of all points
  where forces are applied and that avoids $p$ obstacles, $p_0$ of
  which intersect the boundary of this convex hull,
  can be replaced by a strut net, supporting the same forces, that has at most $q+p-p_0$ elementary loops. 
\end{Thm}

This generalizes the result in  \cite{Bouchitte:2019:FCW} that used a
similar argument to treat
the case where $p=p_0=0$ ($q\ne 0$). \cn

\begin{figure}
	\centering
	\includegraphics[width=0.9\textwidth]{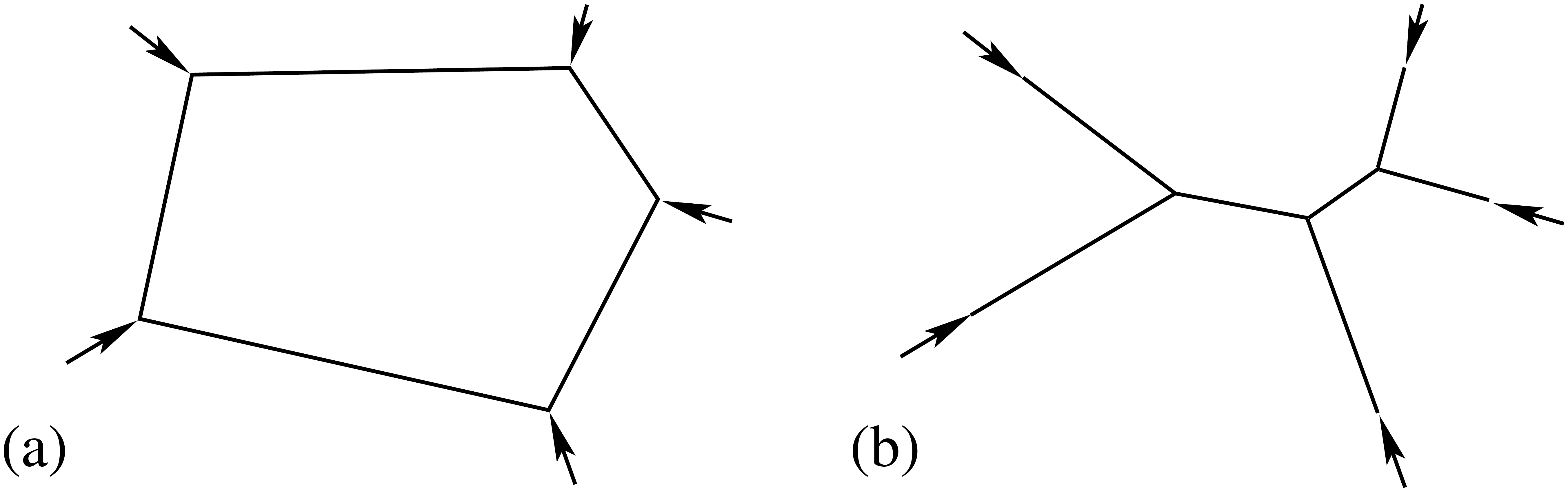}
	\caption{Here (a) is a convex polygonal loop, possibly lying within a strut net that may contain many loops.
          The arrows \crev denote \cn the net forces acting on each vertex. Provided there are no obstacles inside the loop,
          it can be replaced by the open strut net as in (b) thus reducing the number of loops in the whole strut net by one.}
          	\labfig{0}
              \end{figure}

              If a loop contains an obstacle it may be possible to replace that loop by the associated open strut net
              provided that open strut net does
              not collide with the obstacle. \cb In the $q=0$ case \cn it follows that, if a net exists that \crev avoids \cn a single obstacle, then there is a net with none or one loop that also avoids the obstacle.
Specifically, the zero loop case corresponds to the open strut net, and if there is one loop \cb that cannot be reduced, \cn then that loop contains the obstacle. 
Note that, as the struts in the loop are under compression, the loop cannot extend outside
the convex hull of the points where the forces are applied. Thus, if an obstacle intersects the open strut net and extends
outside the convex hull of the points where the forces are applied, then there is no strut net that avoids the obstacle: see
\fig{2}.

\begin{figure}
	\centering
	\includegraphics[width=0.9\textwidth]{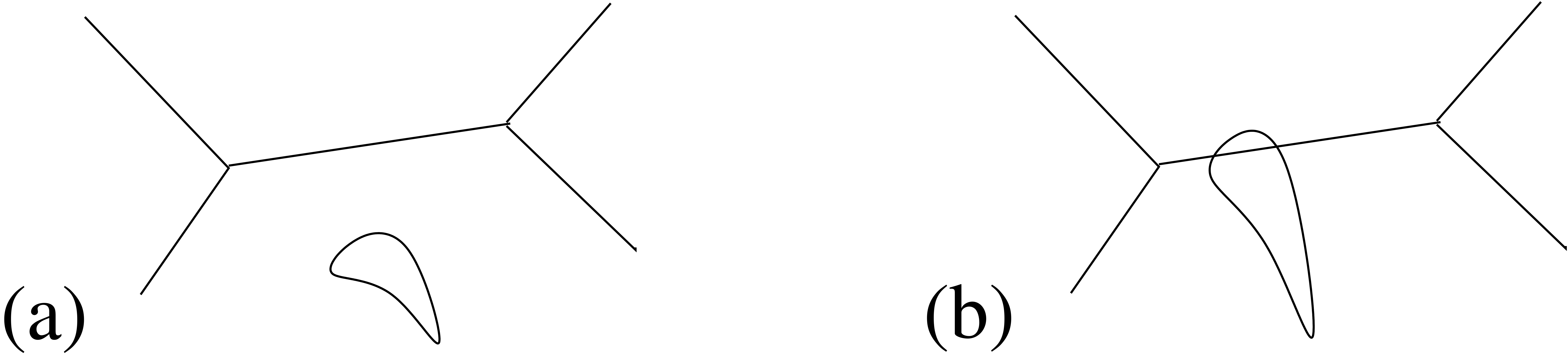}
	\caption{An obstacle may lie partly outside the convex hull of where the forces are applied, but outside the open strut net, as in (a). However if
          the obstacle lies partly outside the convex hull of where the forces are applied and intersects the open strut net, as in (b), then there is
          no strut net that avoids the obstacle.}
	\labfig{2}
      \end{figure}

\section{Solution for the obstacle problem with forces at the boundary of a convex polygon} \labsect{3}
\setcounter{equation}{0}

Our objective is to find an algorithm for determining if a strut net exists that avoids a given set of simply connected obstacles and which
supports a given set of forces at the boundary of a convex polygon. \crev  We will first assume that there is just one obstacle and we will enlarge the region it can occupy. \cn Here by enlargement we use
the following definitions:

\begin{Def}
  Given a strut net \cb avoiding \cn an obstacle $\CO$, the region available
  to $\CO$ is the interior of the strut loop, having no internal
  struts, that contains $\CO$.
  \label{def2}
\end{Def}

\begin{Def}
  Given a region $\GG$ \cb available to $\CO$ for some strut net, then $\GG'$ is an enlargement of it if  $\GG'$ is available to $\CO$ for some strut net and 
  $\GG'\supseteq\GG$. \cn
  \label{def1}
\end{Def}

Thus, if the obstacle is in a region $\GG$, then it is surely in the enlarged region $\GG'$: $\CO\subseteq\GG$ implies $\CO\subseteq\GG'$.
While the procedure will be illustrated by figures in which there are only 4 applied forces, the analysis holds when
there are any number of forces applied at the vertices of a convex polygon.

Given a strut net, and associated strut net function $\psi(\Bx)$, we can modify the strut net to try to accommodate, or better accommodate,
an obstacle $\CO$ by replacing $\psi(\Bx)$ with
\beq \phi(\Bx)=\min\{\psi, L\},\text{  with  }L(\Bx_k)\geq L_k(\Bx_k)\text{ for all }k, \eeq{0.10}
where $L(\Bx)$ is a linear function to be chosen so that $\phi(\Bx)=L(\Bx)$ on $\CO$. Thus, $\psi(\Bx)\geq L(\Bx)$ on $\CO$. 
We call $L(\Bx)$ the cleaving plane. Now, if for any $\Bx$,
\beq L(\Bx)\leq \psi(\Bx), \eeq{0.11}
then surely because $\psi(\Bx)\leq\phi_0(\Bx)$, one has that
\beq L(\Bx)\leq \phi_0(\Bx). \eeq{0.12}
Thus, the region $\GG$ available to the obstacle is enlarged if the
cleaving plane cleaves the open \crev strut net function \cn
rather than any other strut net function
supporting the forces: see \fig{1}(a)-(d).
\begin{figure}
	\centering
	\includegraphics[width=0.9\textwidth]{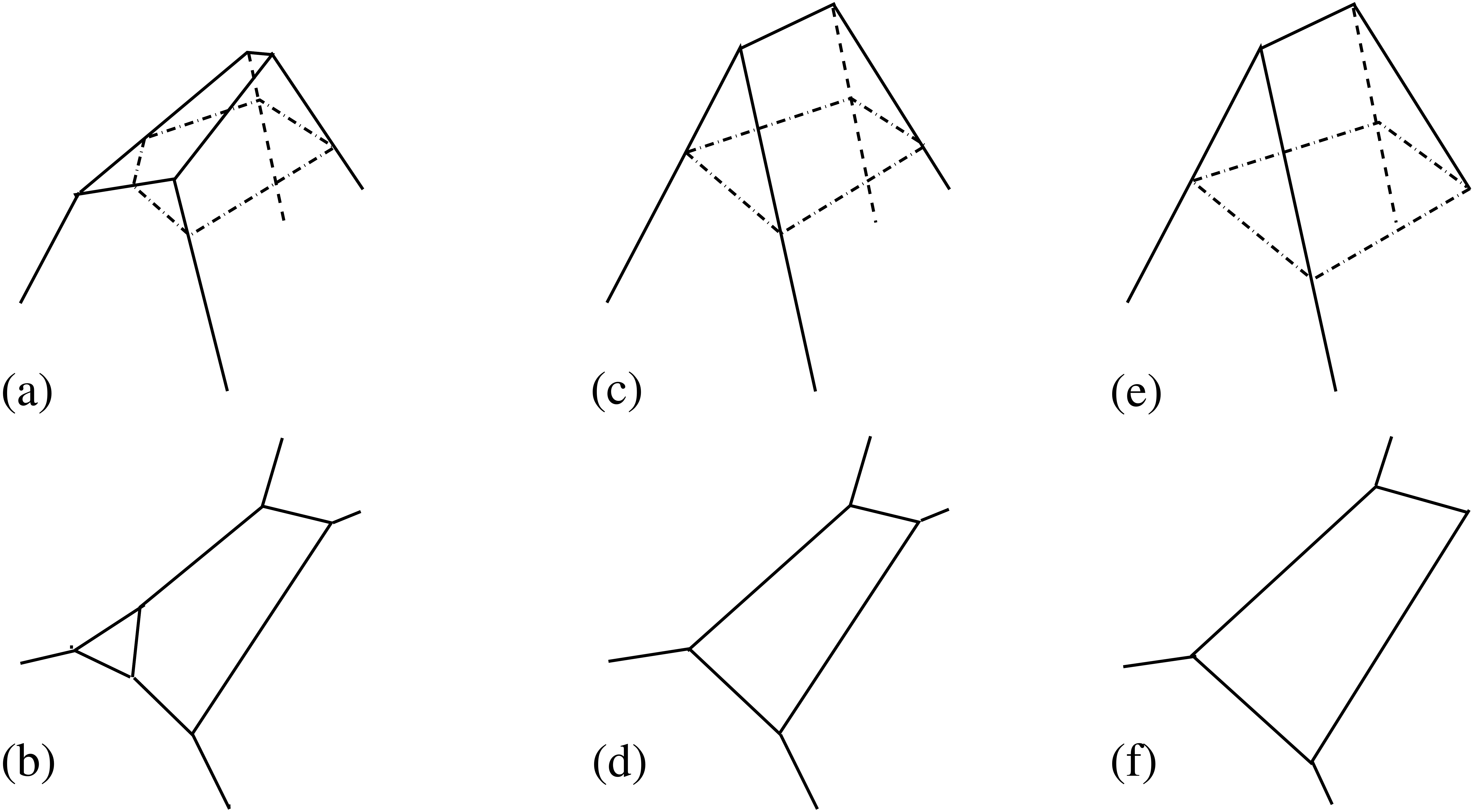}
	\caption{Here (a) shows the cleaving (dot-dash line) of a strut net function (solid and dashed lines), resulting in the strut net (b). If
          instead as in (c) one cleaves, with the same cleaving plane, the open strut net function, one enlarges the region that the
          obstacle can occupy, and the corresponding strut net (d) has a single \crev elementary internal loop. \cn Finally, by moving the cleaving plane downwards keeping the same orientation
          until it touches one of points $L_k(\Bx_k)$, as in (e), one further enlarges the region that the obstacle can occupy, as in (f).}
	\labfig{1}
      \end{figure}

\crev Hence,  the existence of a strut net avoiding $\CO$ is equivalent to the existence of a linear function $L(\Bx)$ satisfying
\beq L (\Bx_i)\geq L_i(\Bx_i) \quad\text{and }\quad L(\Bx)\leq L_i (\Bx)\text{ on } \CO \text{ for all }i. \eeq{0.13}
\cn
\subsection{One or more obstacles which are convex polygons, or which can be approximated by convex polygons}
Note that, \crev for any cleaving plane $L(\Bx)$, the set $\{\Bx:\phi_0(\Bx)\geq L(\Bx)\}$ is a convex set because of the concavity of the open strut net function $\phi_0(\Bx)$. \cn 
Hence, there is no loss of generality by assuming that \crev $\CO$ \cn is convex. Therefore, consider the case when $\CO$ is a convex polygon or is approximated by a convex
polygon containing $\CO$. Let $\By_p, p=1,2,\ldots,m$ be the vertices of this polygon. Then, \eq{0.13} will hold if $L(\Bx)$ satisfies
\beq L(\Bx_i)\geq a_i\equiv \Gf_0(\Bx_i) \text{ and } L(\By_p)\leq b_p\equiv\Gf_0(\By_p)\text{ for all }i,p.
\eeq{z.1}
Now, $L(\Bx)$ can be written as $L(\Bx)=\Bv\cdot \Bx+c$ and we are
looking for $\Bv$ and $c$ \crev so that the inequalities \eq{z.1} are satisfied, that is, for all $i$ and $p$, \cn
\beq \Bv\cdot \Bx_i +c \geq a_i,\quad \Bv\cdot \By_p+c \leq b_p. \eeq{z.2}
The question becomes whether there is a feasible \cb domain \cn  of $(\Bv, c)$ pairs satisfying
this system of linear inequalities: this is a standard problem in linear programming theory.

Note that this is easily generalized to the case of multiple inclusions $\CO_q,\, q=1,2,\ldots, s$, with associated cleaving planes
$L^{(q)}(\Bx)$. Then, \eq{0.13} is replaced with
\beq L^{(q)} \geq L_i(\Bx_i) \quad\text{and }\quad L^{(q)}\leq L_i,\quad L^{(q)}\leq L^{(r)}\text{ on } \CO_q, \text{ for all }i,q,r.
\eeq{z.3}
Treating each  $\CO_q$ as convex or approximating it by a convex
polygon containing $\CO_q$, we let $\By_p^{(q)}, p=1,2,\ldots,m(q)$ be its vertices. Then \eq{z.3} will hold if the cleaving
planes $L^{(q)}(\Bx)$ satisfy for all $i,p,r$
\beq L^{(q)}(\Bx_i)\geq a_i \equiv \Gf_0(\Bx_i)=L_i(\Bx_i), \eeq{z.3a}
and
\beq
L^{(q)}(\By_p^{(q)})\leq b_p^{(q)}\equiv\Gf_0(\By_p^{(q)}), \quad
L^{(q)}(\By_p^{(q)})\leq L^{(r)}(\By_p^{(q)}).
\eeq{z.4}
Writing $L^{(q)}(\Bx)$ as $L^{(q)}(\Bx)=\Bv^{(q)}\cdot \Bx+c^{(q)}$ we arrive at the linear system of inequalities
\beq  \Bv^{(q)}\cdot \Bx_i +c^{(q)} \geq a_i,\quad \Bv^{(q)}\cdot \By_p^{(q)}+c^{(q)} \leq b_p^{(q)}, \quad  \Bv^{(r)}\cdot \By_p^{(q)}+c^{(r)}-\Bv^{(q)}\cdot \By_p^{(q)}+c^{(q)} \geq 0,
\eeq{z.5}
which must hold for all $i,q$ and $r$. This is again a standard problem of determining if there is a non-empty feasible \cb domain \cn  in the $3s$ dimensional space consisting
of $s$ pairs $(\Bv^{(q)}, c^{(q)})\, q=1,2,\ldots, s$. If it is non-empty,
then associated with a  \cb $(\Bv^{(q)}, c^{(q)})\, q=1,2,\ldots, s$ \cn in the feasible \cb domain \cn  are cleaving planes $L^{(q)}(\Bx),\, q=1,2,\ldots, s$ and when $\phi_0(\Bx)$ is cleaved by them we obtain the strut net avoiding the obstacles.  

\subsection{Enlarging as much as possible the region available to the obstacle} 

Let us return back to the case where there is a single obstacle and assume that the feasible \cb domain \cn  associated with \eq{z.2} is non-empty. While the feasible \cb domain \cn 
allows us to identify all strut nets with one
strut loop that contains $\CO$, it does not identify those strut nets with as much room as possible around $\CO$. By this we mean that the region $\GG$ available to $\CO$
is \cb maximal, \cn in the sense that if $\GG'$ is the region available to $\CO$ in another strut net and $\GG\subseteq\GG'$,
then $\GG'=\GG$. Note that $\GG$ is not necessarily unique, and a different $\GG$ may have less room around the obstacle on one side and
more on another side. Here we will show that if $\CO$ lies in a \cb maximal \cn $\GG$ then
\beq  L(\Bx_k)=L_k(\Bx_k)\text{ and }L(\Bx_\ell)=L_\ell(\Bx_\ell), \eeq{0.13aa}
for some $k$ and $\ell$.

If, for any $\Bx$, \eq{0.12} holds then surely
\beq L'(\Bx)\leq \phi_0(\Bx)\text{  with }  L'=L-\min_i (L(\Bx_i)-a_i),\quad a_i=L_i(\Bx_i), \eeq{0.14}
because $L'(\Bx)\leq L(\Bx)$.
In other words, the cleaving plane can be lowered while keeping its orientation fixed, expanding the region $\GG$ that can be occupied by the
obstacle. This can be continued until this cleaving plane $L'$ first touches a point $A=L_k(\Bx_k)$, where $i=k$ attains the
equality in \eq{0.14}: see \fig{2}.

Now we tilt the cleaving plane about the line through $A$ and another point $B\in L\cap L_k$  with $B\ne A$. This line is where $L'(\Bx)$ and $L_k(\Bx)$ intersect.
This gives a new cleaving plane
\beq L''=L'+\alpha(L'-L_k), \text{ with }\Ga>0. \eeq{0.15}
 We can continue this tilting by increasing $\Ga$ until
\beq \alpha=\min_{j} \frac{L'(\Bx_j)-L_j(\Bx_j)}{L_k(\Bx_j)-L'(\Bx_j)}, \eeq{0.16}
when $L''(\Bx)$ first touches another point $C=L_\ell(\Bx_\ell)$ in which $j=\ell$ achieves the equality in  \eq{0.16}.
It also touches $A$. Note that if for any $\Bx$ the first inequality in \eq{0.14} holds then surely
\beq L''(\Bx)\leq \phi_0(\Bx), \eeq{0.17}
because that first inequality implies $L'(\Bx)\leq L_k(\Bx)$ at that point. This implies that as the \cred cleaving plane is tilted
until it touches $B$, the region $\GG$ that can be occupied by the obstacle enlarges further: see \fig{3}(a)-(d).

One can subsequently rotate the cleaving plane about the line joining $A$ with $L_\ell(\Bx_\ell)$ until the cleaving plane meets
$L_h(\Bx_h)$ for some $h$. In our example, rotation one way gives \fig{3}(e) with $L_h(\Bx_h)=D$
and the corresponding strut net \fig{3}(f), while rotation the other way gives \fig{3}(g) with $L_h(\Bx_h)=E$
and the corresponding strut net \fig{3}(h). However, this rotation does
not generally enlarge the region $\GG$ that can be occupied by the obstacle.
The rotation tips the cleaving plane up on one side of the line $L_k(\Bx_k)$ with $L_\ell(\Bx_\ell)$ and
down on the other side. The region $\GG$, with the strut loop as its boundary, shrinks on the up side and expands on the down side. The exception,
as illustrated in \fig{3}(i) with the corresponding strut net \fig{3}(j) is when $\ell=k\pm 1$, i.e. when $L_\ell(\Bx_\ell)$ is
on the same facet of $\phi_0(\Bx)$ as $L_k(\Bx_k)$. Then the region $\GG$ lies on
one side of the line joining $\Bx_k$ with $\Bx_\ell$ and, as illustrated in \fig{3}(k) we may rotate the cleaving plane about the line so the plane
tilts down on the side where $\GG$ is located and touches $D=L_h(\Bx_h)$, producing the strut net of \fig{3}(l).

\begin{figure}
	\centering
	\includegraphics[width=0.9\textwidth]{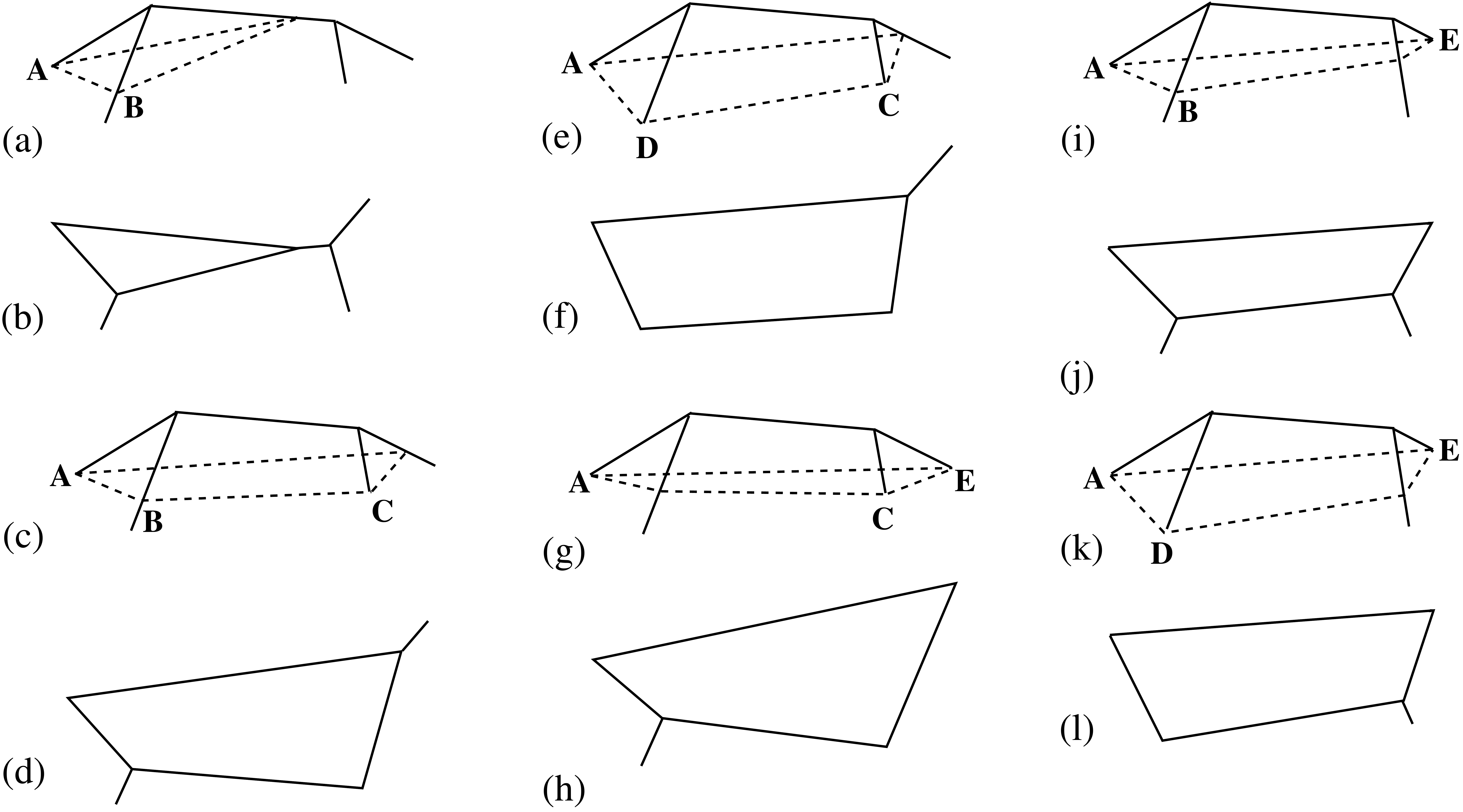}
	\caption{Various configurations, (a), (c), (e), (g), (i), and (k) of the cleaving plane, denoted by dashed lines,
          that cleave $\phi_0(\Bx)$, denoted by solid lines. The associated strut nets are below each of these subfigures. The procedure of tilting
          and then rotating the cleaving plane follows that described in the text. 
        }
	\labfig{3}
      \end{figure}

      This analysis provides an algorithm to generate all \cb maximal \cn regions $\GG$. One starts with each $i=1,2,\ldots,n$ and
      takes the cleaving plane through $L_i(\Bx_i)$, $L_{i+1}(\Bx_{i+1})$ and a third point
      $L_j(\Bx_j)$ (with no points  $L_h(\Bx_h), h=1,2,\ldots,n$ lying above the plane). One next rolls the cleaving plane,
      rotating it about the line joining  $L_i(\Bx_i)$ and $L_j(\Bx_j)$ until it first hits another point $L_k(\Bx_k)$.
      Then one rolls it further about the line joining  $L_i(\Bx_i)$ and $L_k(\Bx_k)$ until it first hits another point
      $L_\ell(\Bx_k)$. This is continued until the rolling cleaving plane hits  $L_{i-1}(\Bx_{i-1})$. As the rolling progresses
      the region enclosed by the strut net loop created by the cleaving plane is \cb maximal \cn and the region available to
      an obstacle placed within it cannot be enlarged. 

\section{Including reactive forces with multiple obstacles} \labsect{4}
\setcounter{equation}{0}

Let us suppose that a subset ${\mathcal I}\subset \{1,\ldots,n\}$ of the $n$ forces are reactive, while the remaining ones (corresponding to the subset of indices ${\mathcal K}:= \{1,\ldots,n\}\setminus {\mathcal I}$) are given.  That means that the forces are not fixed at the points $\Bx_{i}$ for $i\in {\mathcal I}$ but that any forces at these points  are admissible as soon as they balance through a strut net the forces $\Bt_{i}$ applied at points $\Bx_{i}$ for $i\in {\mathcal K}$ .
The planes $L_k$ are no longer fixed by equations \eq{A.1}. 
When $i\in {\mathcal I}$, these equations, determine the reactive force $\Bt_i$:
\beq \Bt_i\equiv -\BR_\perp(\Grad L_{i+1}-\Grad L_{i}),\quad i\in {\mathcal I} \eeq{4.2a}
and provide no constraints on the unknown $L_{i+1}$ and $L_{i}$
while, for the remaining $i\in {\mathcal K}$, they provide the constraints
\beq \BR_\perp(\Grad L_{i+1}-\Grad L_{i})=-\Bt_i, \quad i\in {\mathcal K}.  \eeq{4.3a}
Since we can add any linear function to the Airy stress function, we are free to
assume, as before, that 
\beq L_1=0.\eeq{4.3aa}
Continuity of the Airy function requires that, for all $i\in \{1,2,\ldots, n\}$,
\beq L_{i+1}(\Bx_i)=L_{i}(\Bx_i).
\eeq{4.3b}
Its concavity requires that, for all $i\not=j$ in $\{1,2,\ldots, n\}$,
\beq L_{j}(\Bx_i)\geq L_{i}(\Bx_i).
\eeq{4.3bb}
As in \eq{z.5}, the existence of a strut net avoiding the  $s$ obstacles $\CO_1$, $\CO_2$, \dots $\CO_s$, each $\CO_q$ being the convex hull of points \crevv $(\By^q_1,\By^q_2,\dots,\By^q_{n(q)})$, \cn
is equivalent to the existence, for all $q\in\{0,1,\dots, s\}$, of linear functions $L^{(1)}$, $L^{(2)}$, \dots $L^{(s)}$ satisfying, for any $i$ in $\{0,1,\dots, n\}$, any $q\not=r$ in $\{0,1,\dots, s\}$ and any $p$ in $\{0,1,\dots, n(q)\}$,
\beq  L^{(q)}(\Bx_i)  \geq  L_i(\Bx_i),
\quad L_{i}(\By_p^{(q)}) \geq L^{(q)}(\By_p^{(q)}), 
\quad L^{(r)}(\By_p^{(q)}) \geq L^{(q)}(\By_p^{(q)}). \eeq{4.3c}

\cb These results, which also generalize the results obtained in Section \sect{3}, can be summarized in the following theorem:
\begin{Thm} The existence of  a strut net balancing forces
$\Bt_{i}$ applied at points $(\Bx_{i})_{i\in {\mathcal K}}$ with the help of reacting points $(\Bx_{i})_{i\in {\mathcal I}}$ and avoiding obstacles $\CO_1$, $\CO_2$, \dots $\CO_s$, each $\CO_q$ being the convex hull of points $(\By^q_1,\By^q_2,\dots,\By^q_{n(q)})$, is equivalent to the existence of linear functions $(L_i)_{i=1}^n$ and $(L^{(q)})_{q=1}^s$ satisfying, for any $i$ in $\{0,1,\dots, n\}$, any $q\not=r$ in $\{0,1,\dots, s\}$ and any $p$ in $\{0,1,\dots, n(q)\}$, the constraints \eq{4.3a}, \eq{4.3aa}, \eq{4.3b} , \eq{4.3bb} and \eq{4.3c}.
Setting
\beq L_i(\Bx)=\Bw_i\cdot \Bx+d_i\quad \text{and}\quad L^{(q)}(\Bx)=\Bv^{(q)}\cdot \Bx+c^{(q)}, \eeq{4.1aa}
the existence of such a strut net is equivalent to the existence of a solution $(\Bw_i,d_i)_{i=1}^n$ and $(\Bv^{(q)},c^{(q)})_{q=1}^s$ of the following linear programming problem
\begin{align}
&  \BR_\perp(\Bw_{i+1}-\Bw_{i})=-\Bt_i,
& \text{for $i\in {\mathcal K}$,} \nonum
&\Bw_1=0,\quad d_1=0,  \nonum
&\Bw_{i+1}\cdot \Bx_i+d_{i+1}=\Bw_{i}\cdot \Bx_i+d_i,
&\text{for $1\leq i\leq n$,} \nonum
&\Bw_{j}\cdot \Bx_i+d_{j}\geq \Bw_{i}\cdot \Bx_i+d_i,
&\text{for $i\not= j$ in $\{1,\dots,n\}$,} \nonum
& \Bv^{(q)}\cdot \Bx_i +c^{(q)}  \geq  \Bw_{i}\cdot \Bx_i+d_i,
&\text{for $1\leq i\leq n$, $1\leq q\leq s$,} \nonum
&\Bw_{i}\cdot \By_p^{(q)}+d_i \geq \Bv^{(q)}\cdot \By_p^{(q)}+c^{(q)},
&\text{for $1\leq i\leq n$, $1\leq q\leq s$, $1\leq p\leq n(q)$,} \nonum
&\Bv^{(r)}\cdot \By_p^{(q)}+c^{(r)} \geq \Bv^{(q)}\cdot \By_p^{(q)}+c^{(q)},
  &\text{for $q\not= r$ in $\{1,\dots,s\}$, $1\leq p\leq n(q)$.}
    \label{LP}
\end{align}
\end{Thm}

Recall that, as soon as admissible quantities $(\Bw_i,d_i)_{i=1}^n$ have been chosen, the reacting forces are determined by  \eq{4.2a}.

We must emphasize the fact that, when a solution exists for the previous linear programming problem, generally an infinity of solutions exists.
Indeed we have not fixed any objective function to minimize. We can fix in very different ways such a linear objective function.
But even when  it is given uniqueness is not guaranteed. This fact is clearly illustrated by the following example. Assume that all bars in the strut net
have a \cg thickness \cn proportional to the force they carry, in such a way that the stress remains constant in the whole structure. Then the total volume $\mathcal V$ of the structure
is proportional to the integral over $\Omega$ of the Laplacian of the strut net function (a negative measure). This linear function of the unknowns, \cg which from here onwards
we will call the total weight of the structure, \cn seems to be a good
candidate for the objective function. However one can notice that, using the divergence theorem, $\mathcal V$ can be computed from the normal derivative of the
strut net function on the boundary of the domain. Therefore it does not involve the $L^{(q)}$ functions: its minimization can help fixing the $L_i$ functions,
that is the reactive forces, but it cannot entirely determine a unique strut net.

We can also treat the case when a subset $\mathcal J \subseteq \mathcal I$ of the points at which the reactive forces act, rather than being fixed,
are confined to some line segments. These segments, for example, could be sections of supporting walls or ground. Then, there are the constraints
\beq \Bx_j\cdot(\BR_\perp\Bg_j)=z_j,\quad  g^-_j\leq \Bx_j\cdot\Bg_j \leq g^+_j, \text{ for } j\in \mathcal J, \eeq{4.1}
where the unit vectors $\Bg_j$ and constants $z_j$, $g^-_j$ and $g^+_j$ are given, but only defined for $j\in{\mathcal J}$.
These segments should be such that the $\Bx_i, i=1,2,\ldots, n$ are still the vertices, going anticlockwise, of a convex polygon for all
choices of points $\Bx_j, j\in{\mathcal J}$ anywhere on the line segments. 
Now we need to add the $\Bx_j, j\in{\mathcal J}$ to the unknowns. Then, those equations in \eq{LP} that involve  $\Bx_j, j\in{\mathcal J}$
provide quadratic (and non-convex) constraints on the unknowns,
and we are left with determining the feasible \cb domain \cn  associated with a quadratic programming problem, involving \eq{LP} and \eq{4.1}, again a standard problem.

It could be the case that reactive forces act at, say, a pair of neighboring points
$\Bx_{i}$ and $\Bx_{i+1}$ with $i\in{\mathcal J}$ and $i+1\in{\mathcal J}$ and share the same line segment. Then,
\beq \Bg_i=\Bg_{i+1}, \quad z_i=z_{i+1}, \quad  g^-_i=g^-_{i+1}, \quad g^+_i=g^+_{i+1},
\eeq{4.7}
and to maintain the anticlockwise order of points \eq{4.1} needs to be
supplemented by the additional constraint that
\beq  (\Bx_i-\Bx_{i+1})\cdot\Bg_i\geq 0, \eeq{4.8}
by, if necessary, reversing the direction of $\Bg_i$, the signs of $z_i$, $g^-_i$ and $g^+_i$, and reversing the inequalities in \eq{4.1},
while maintaining \eq{4.7}.
The generalization to the case of more reactive forces acting at neighboring points sharing the same interval is
straightforward.

There is an alternative to allowing the points $\Bx_{j}, j\in{\mathcal J}$ to range along these line segments: one can distribute along each line segments sufficiently many additional
fixed points $\Bx_i$ at which reactive forces act. This has the advantage of replacing the quadratic programming problem with a linear one, generally at the sacrifice of having more unknowns. 
This is what has been done in the numerical examples in Sect. \sect{Subsec_Construction_strut_nets} where 200 points with reactive forces have been distributed along two supporting segments.


\section{Numerical Results} \labsect{5}
     \setcounter{equation}{0}

\crev The first part of \cn\cb this  \cn\crev section deals with some numerical  applications  of Theorem 1 of Section \sect{2}, the second part treats the issue of enlarging the region where the obstacle can be placed, while the last part discusses numerical examples of  the
linear programming procedure described in \cn\cb Theorem 2 \cn\crev of Section \sect{4}. Specifically, in Section \sect{Sec_Reduction_numerics}, we consider two examples in which the internal loops of a strut net that avoids an obstacle are simplified into open \cn \cb strut nets, \cn \crev up to a certain number of irreducible elementary loops. Indeed, according to Theorem 1 in Section \sect{2}, the number of loops that cannot \cn \cb further be simplified depends on the number of points of application of the forces that lie inside the convex hull and on the position of these forces and the obstacle itself. \cn\crev In Section \sect{Sub_sec_enlaring_gamma}, we implement the algorithm,  presented at the end of Section \sect{3}, that allows one to enlarge the area available to place an obstacle. \cn Finally, in Section \sect{Subsec_Construction_strut_nets}, \cn strut nets are generated that support sets of active applied forces at given
points and reactive forces at other given points that act in response
to the applied forces. 
The first example gives a strut net with two triangular cells that
supports a single active force and avoids four obstacles.  Four 
subsequent examples examine strut net models for masonry arches
modeled as rigid no-tension bodies \cite{heyman1995,marmo2017}, which are subject to different
loading conditions and are required to avoid various
obstacles. 
\crev In all numerical examples, we use abstract units for lengths and forces.\cn


\subsection{Loop reduction}\labsect{Sec_Reduction_numerics}
\subsubsection{Reduction to a funicular arch strut net} \labsect{5loop}

\begin{figure}[h!]
	\centering
	\includegraphics[width=1.0\textwidth]{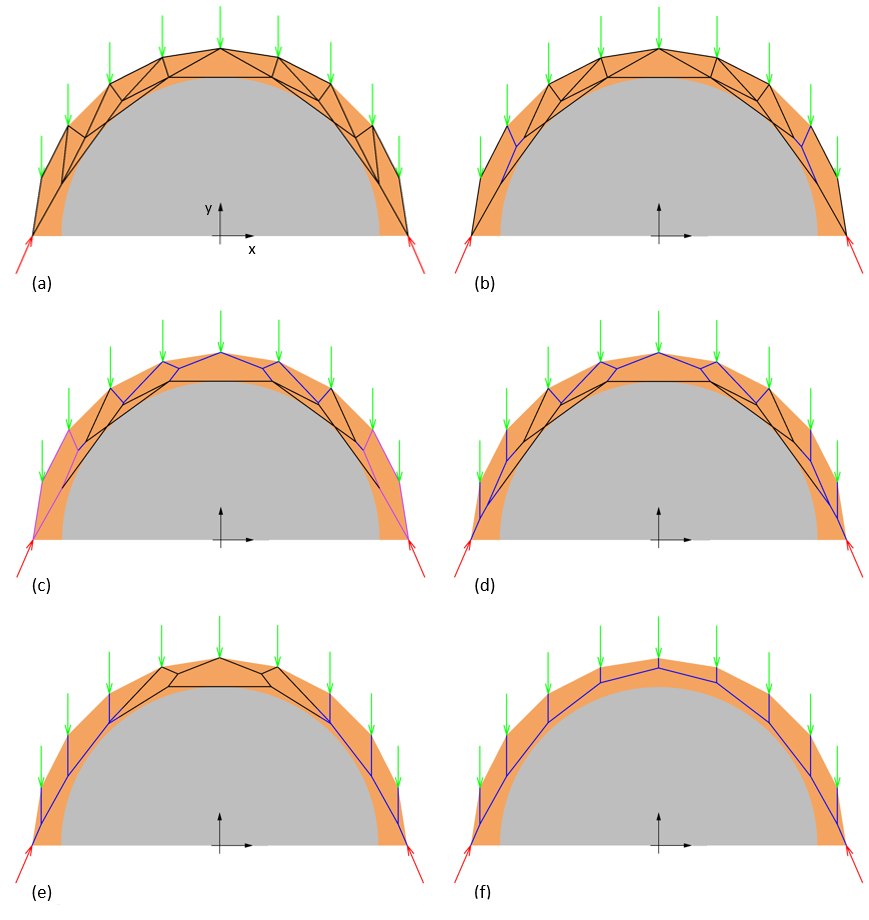}
	\caption{ An arched strut net that avoids a half disk
		(represented in gray) and supports 9 vertical forces of
		equal magnitude acting downwards \cn  along the arch, with the help of two 
		supporting reactive forces (marked in red) acting at the end points of the
		arch. \cn  An open strut net is obtained from an initial net exhibiting  $\ell=13$ closed loops (panel (a)), through the loop reduction procedure presented in \cite{Bouchitte:2019:FCW} (panels (b)-(f)).
		The masonry arch containing the \cn  sequence \cn  of strut
		nets is colored light blue. \cn  
		\cred Note that the reactive forces at the end points of the arch have been scaled by a factor of 0.2. \cn
		See also Movie S1 in the
		supplementary material. \cn}
	\label{arch_reduction}
\end{figure}

\crevv Let us consider strut net models of a semi-circular masonry arch with horizontal  span \cn\cg of the extrados \cn\crevv
equal to 18, a rise at the intrados of 7.6 and a \cn \cg thickness of 1.4.
\cn\crevv
The arch is loaded by a set of 9 active forces with unit magnitude \cn\cg directed downwards, uniformly distributed at equal angles of  $\pi/20$ \cn\crevv along the extrados, and 2 reactive forces acting at the \cn\cg ends of the arch. \cn
An initial strut net featuring $\ell = 13$ closed loops
was obtained by employing  the \cred numerical procedure \cn presented in \cite{Skelton:2014} under the action of the given active forces (a unique material was assigned to all the struts, see Fig. \ref{arch_reduction}(a)).
Such a strut net was next reduced to an open strut net ($\ell=0$),
employing the loop reduction procedure presented in
\cite{Bouchitte:2019:FCW} (see Figs. \ref{arch_reduction}(b)-(e) and Movie S1).
Since the active and reactive external forces of the current example
are applied at the vertices of a convex polygon ($q=0$), \cn  it is easily
observed that $\ell\leq q+p-p_0$, \cn  in agreement with Theorem
1 of of Section \sect{2}.  
It is easy to show that 
the
open strut net coincides with the funicular polygon of the active
forces, which passes through the end points of the arch and has the
initial and final segments parallel to the lines of the reactive
forces \cite{Allen:1998}.
\crev
Note that the open strut net function \(\phi_0(\Bx)\) for this example, plotted in \fig{arch_reduction_Airy}(a) for \(L_1=0\), provides the open \cb strut net \cn in \fig{arch_reduction_Airy}(b), which coincides with the one represented in Figure \ref{arch_reduction}(f), obtained through loop reduction. 

\cn

\begin{figure}[h!]
	\centering
	\includegraphics[width=0.85\textwidth]{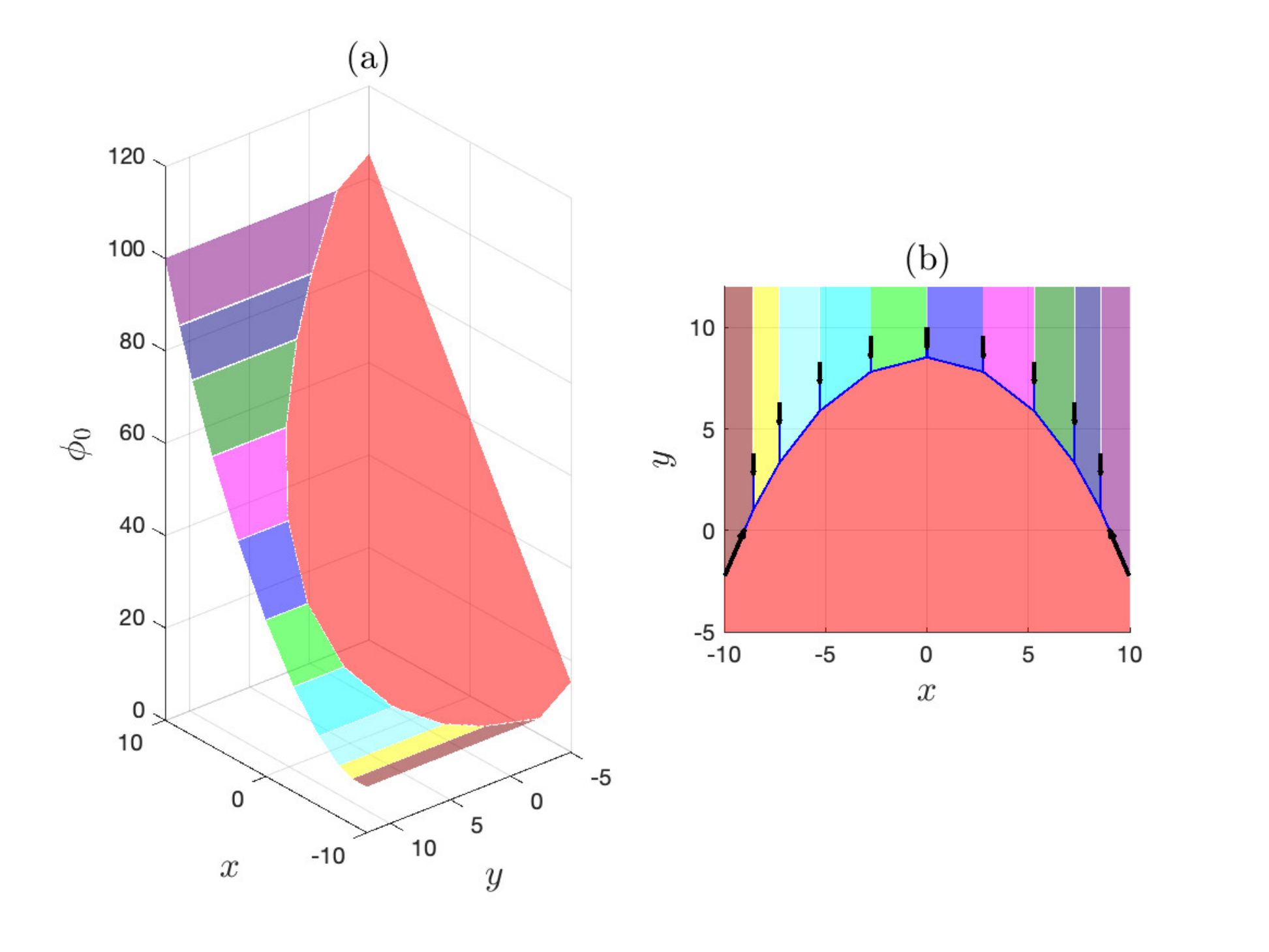}
	\caption{ (a) Here is the open strut net function \(\phi_0(\Bx)\) associated with the example considered in Figure \ref{arch_reduction}. (b) represents the related open
          \cb strut net \cn, that coincides with Figure \ref{arch_reduction}(f). Note that the reactive forces at the end points of the arch have been scaled by a factor of 0.5.}
	\labfig{arch_reduction_Airy}
\end{figure}

\subsubsection{Loop reduction for a strut net with an internal force and that avoids an obstacle} \labsect{Subsec_Reduction_non_convex_obstacle}

Here we consider the example depicted in \fig{loop_reduction}(a), in which one of the forces, \(\Bf_4=(0.874,0.574)\), is applied inside the polygon formed by the remaining forces. Indeed, the point \(\Bx_4=(5,5)\), lies inside the polygon formed by \(\Bx_1=(10,0)\), \(\Bx_2=(0,0)\), \(\Bx_3=(-5,10)\), and \(\Bx_5=(12,7)\) where the forces \(\Bf_1=(2.751,-2.319)\), \(\Bf_2=(-3.415,-1.933)\), \(\Bf_3=(-3.149,1.786)\), and \(\Bf_5=(2.938,1.891)\)  are applied. According to the theory presented in \cite{Bouchitte:2019:FCW}, a \cb strut net \cn supporting such forces with all the elements under compression exists. An example is provided by the strut net connecting all the points pairwise as showed in \fig{loop_reduction}(a). Now suppose that an obstacle is placed inside one of the loops of the \cb strut net \cn. As stated by Theorem 1 in Section \sect{2}, in which $q=1$, $p=1$, and $p_0=0$,  we can replace the \cb strut net \cn connecting the points pairwise by one that has at most $q+p-p_0=2$ elementary loops, as showed in the supplementary Movie S2. The result of the loop reduction is portrayed in \fig{loop_reduction}(b), where the two remaining elementary loops cannot be reduced further.

\begin{figure}[h!]
	\centering
	\begin{subfigure}{.45\textwidth}
	\includegraphics[width=\textwidth]{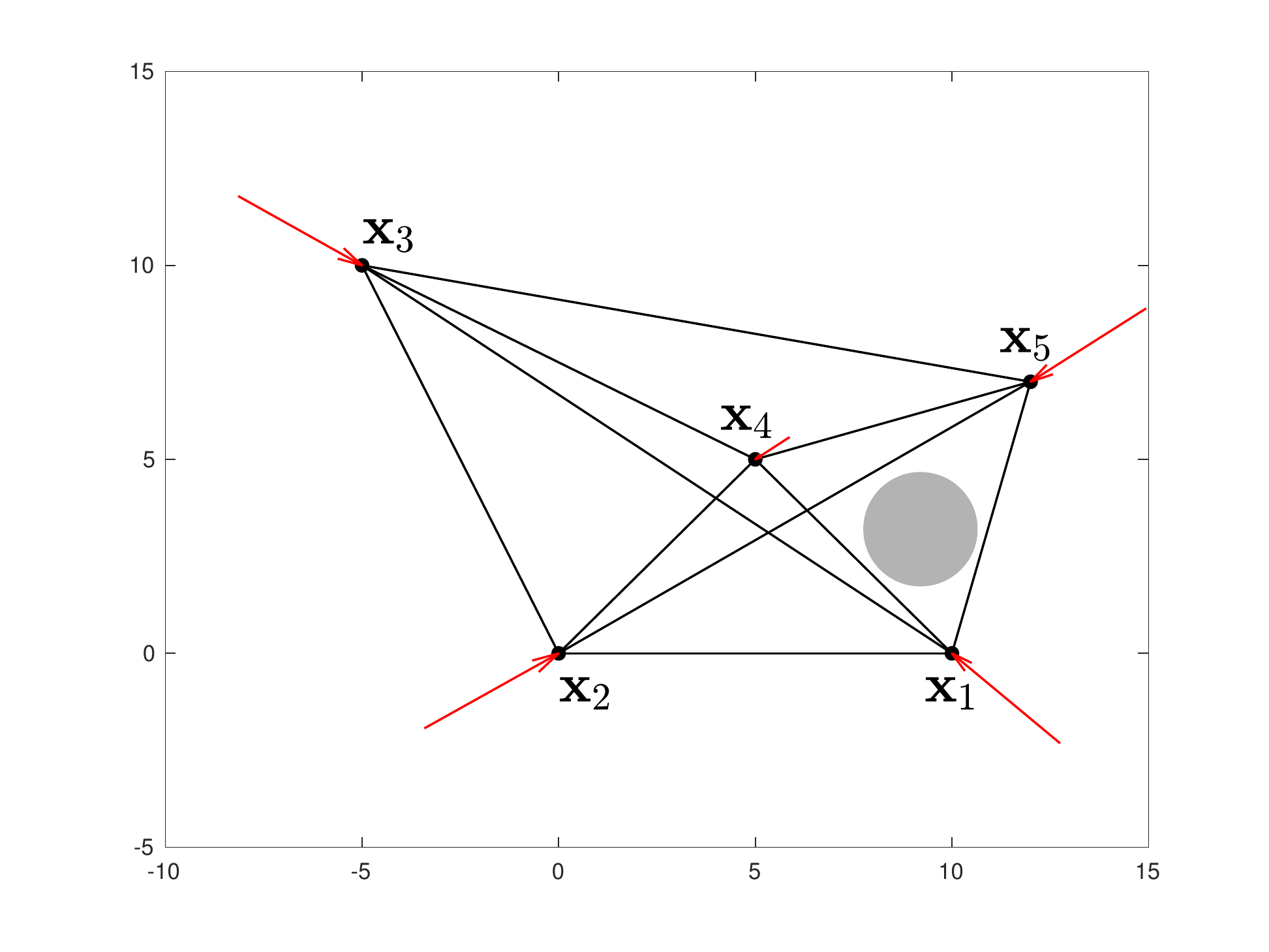}
	\caption{}
	\end{subfigure}
\begin{subfigure}{.45\textwidth}
	\includegraphics[width=\textwidth]{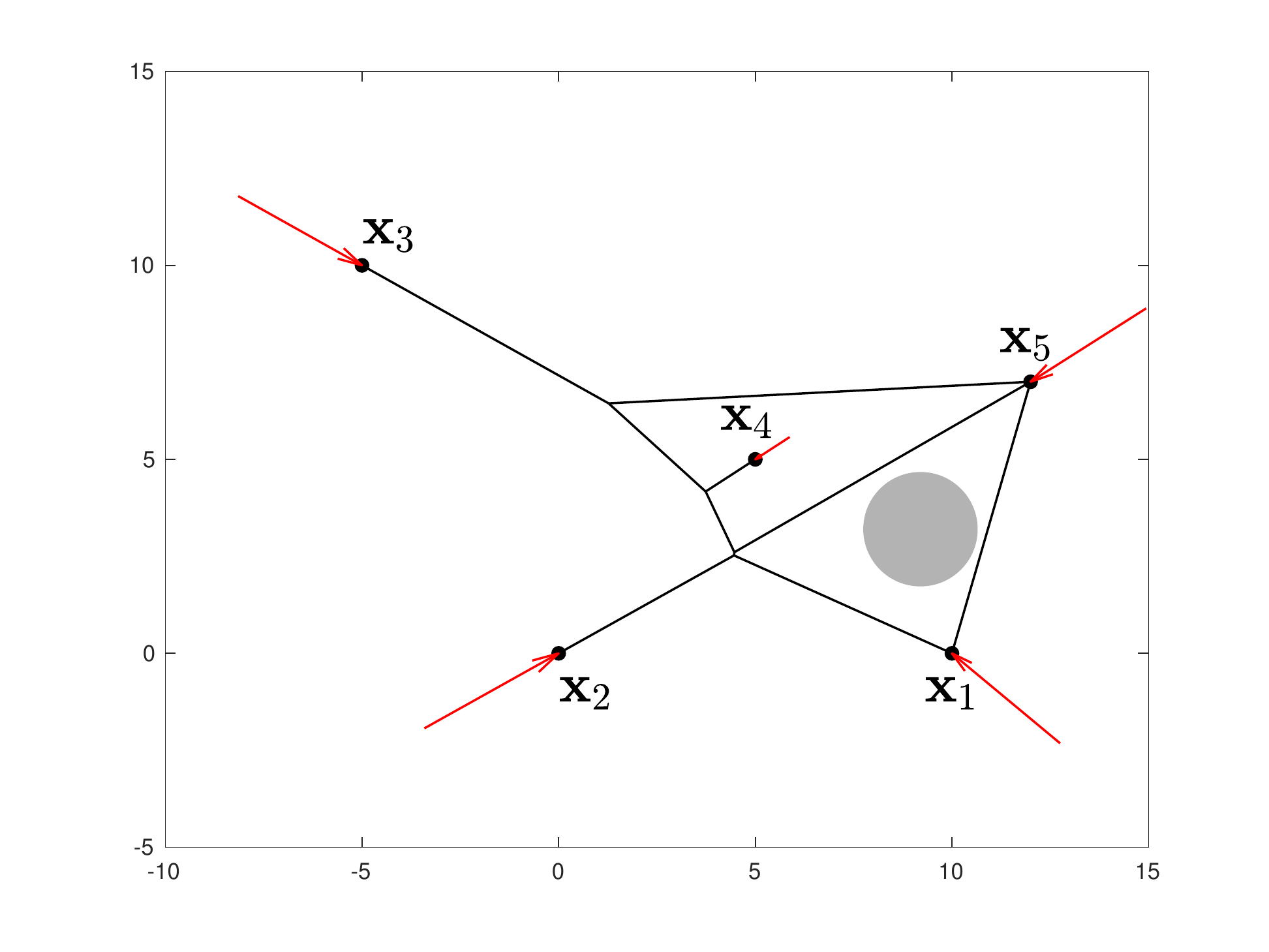}
	\caption{}
\end{subfigure}
	\caption{(a) A \cb strut net \cn that supports the applied forces and avoids the obstacle (here represented by the gray circle) is the one connecting the points pairwise. By applying loop reduction, as illustrated step by step in Movie S2 in the supplementary material, it is possible to reduce the number of elementary loops with at most $q+p-p_0=2$ remaining loops, as illustrated in (b).}
	\labfig{loop_reduction}
\end{figure}

\subsection{Enlarging the region available to the obstacle}\labsect{Sub_sec_enlaring_gamma}

In this section, we provide a numerical example that illustrates the algorithm presented at the end of Section \sect{3} to enlarge the region available to the obstacle. We start with 7 forces, \(\Bf_1=(-1,-4)\), \(\Bf_2=(-3,-1)\), \(\Bf_3=(-2,2)\), \(\Bf_4=(-3,5)\) , \(\Bf_5=(1,1)\) , \(\Bf_6=(6,2)\) ,and \(\Bf_7=(2,-5)\), applied respectively at the points  \(\Bx_1=(0,0)\), \(\Bx_2=(-16,7)\), \(\Bx_3=(-13,16)\), \(\Bx_4=(2,20)\), \(\Bx_5=(12,19)\), \(\Bx_6=(12,13)\), and \(\Bx_7=(10,0)\). The open strut net function for this example is shown in \fig{Airy_7points}(a), and the associated open \cb strut net \cn is represented in \fig{Airy_7points}(b). In order to find the \cb maximal \cn regions \(\GG\), we start with point \(\Bx_1\) and consider the cleaving plane passing through \(L_1(\Bx_1)\), \(L_2(\Bx_2)\), and \(L_3(\Bx_3)\), see \fig{Airy_7points}(c), which provides the region \(\GG\) depicted in black in \fig{Airy_7points}(d). Then, we roll the cleaving plane about the line connecting \(L_1(\Bx_1)\) and \(L_3(\Bx_3)\), until the plane touches the point \(L_4(\Bx_4)\), \fig{Airy_7points}(e)-(f), and this is continued until the plane touches \(L_7(\Bx_7)\), as showed in \fig{Airy_7points}(g)-(h). Movie S3 in the supplementary material shows the rolling of the plane from the initial configuration \fig{Airy_7points}(c) to the final configuration \fig{Airy_7points}(h): as the rolling occurs, the region in which the obstacle can be placed is \cb maximal. \cn 
\begin{figure}[h!]
	\centering
	\begin{subfigure}{0.49\textwidth}
		\centering
		\includegraphics[width=\textwidth]{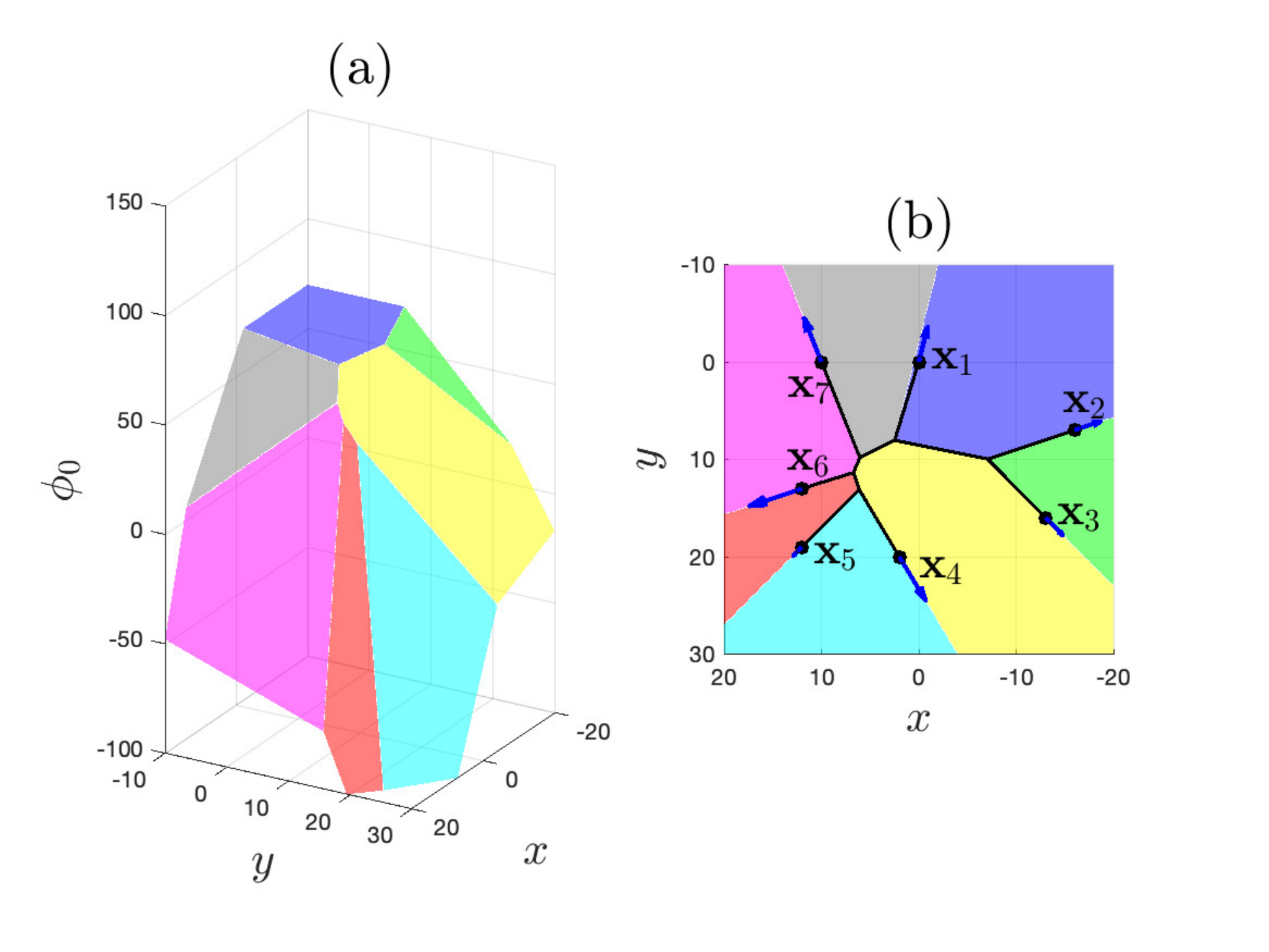}
	\end{subfigure}
	\begin{subfigure}{0.49\textwidth}
		\centering
		\includegraphics[width=\textwidth]{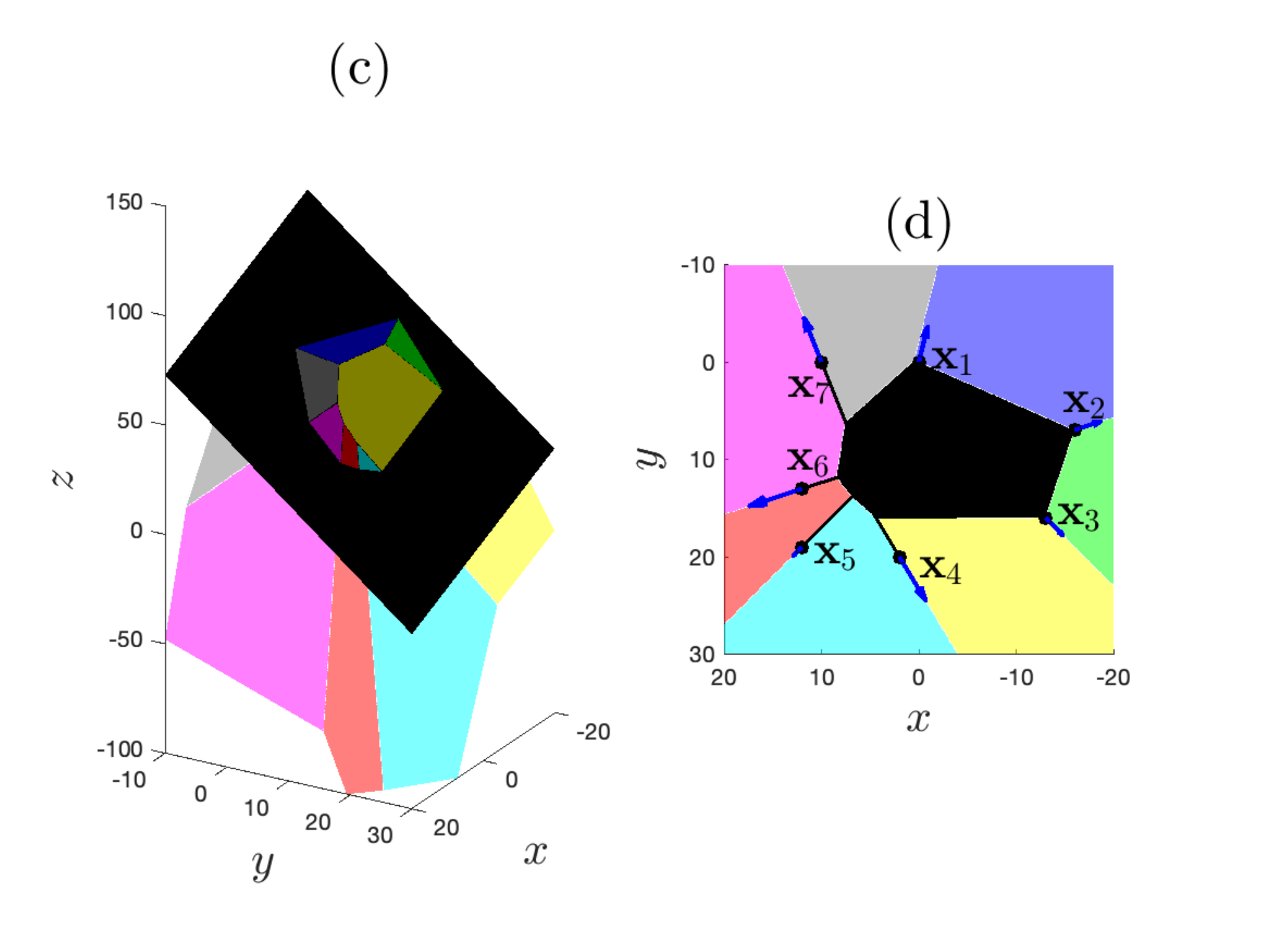}
	\end{subfigure}
\begin{subfigure}{0.49\textwidth}
	\centering
	\includegraphics[width=\textwidth]{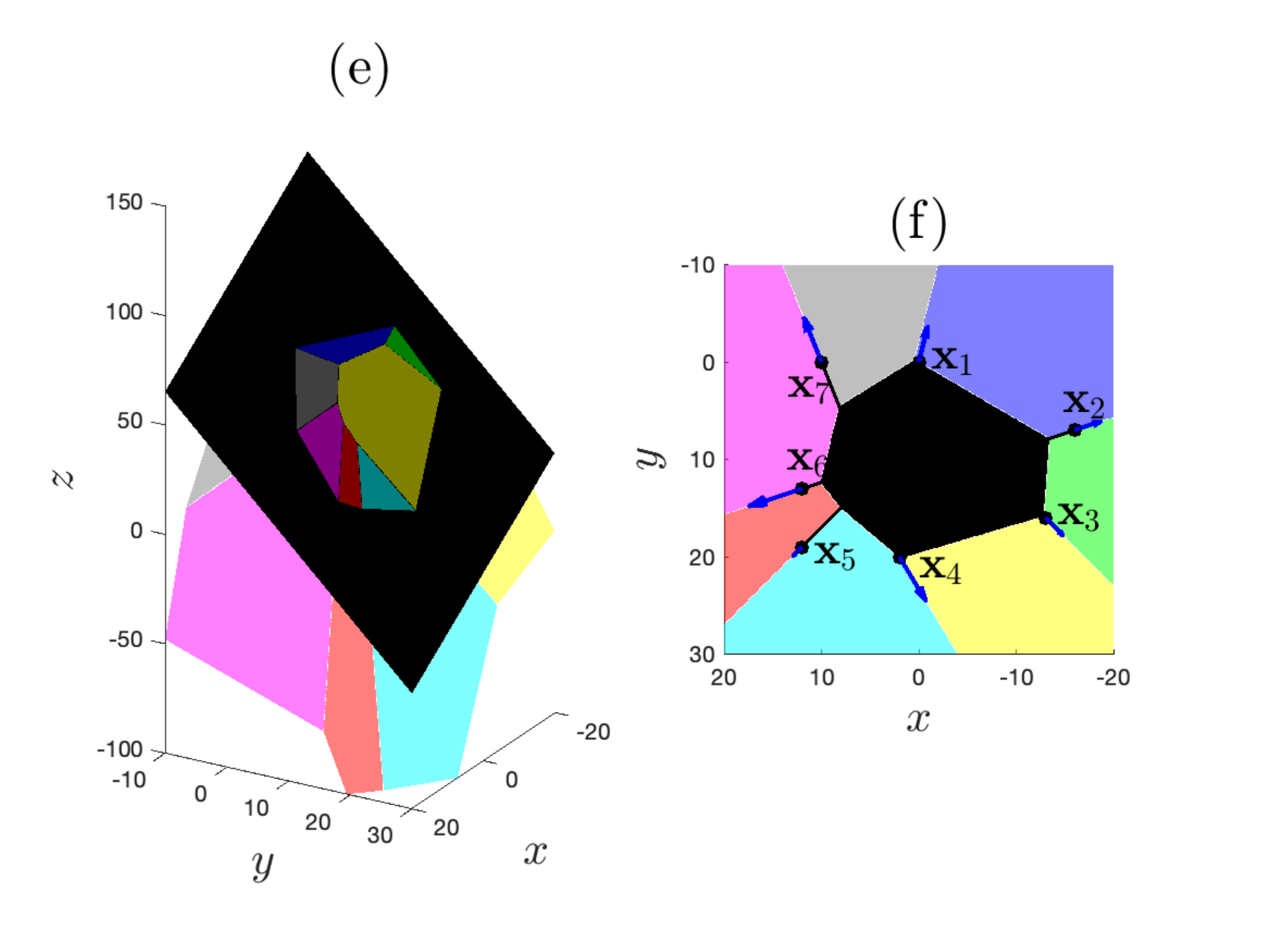}
\end{subfigure}
\begin{subfigure}{0.49\textwidth}
	\centering
	\includegraphics[width=\textwidth]{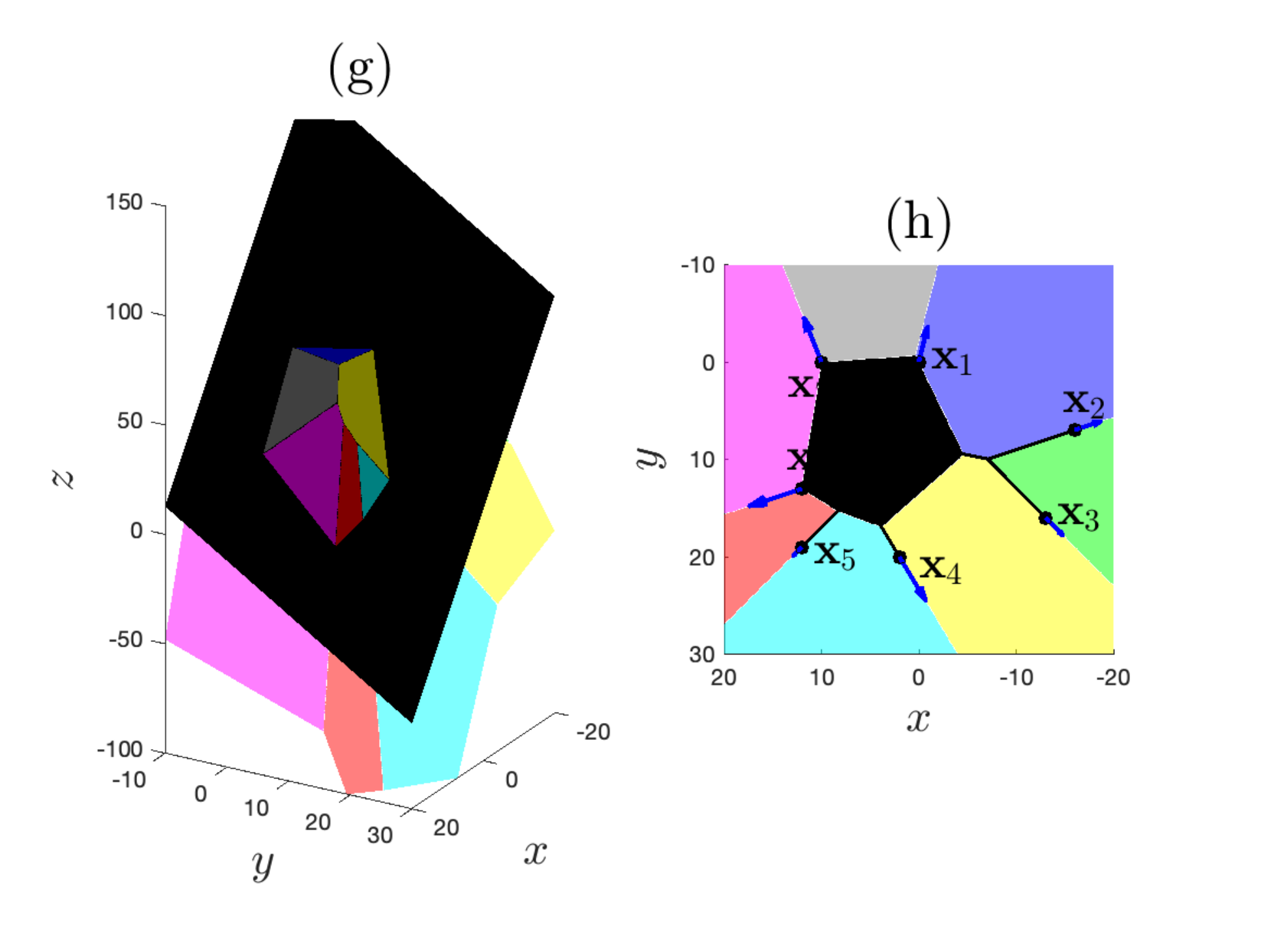}
\end{subfigure}
	\caption{(a) Here's the open strut net function \(\phi_0(\Bx)\) associated with the open \cb strut net \cn depicted in (b), when \(L_1=0\). In (c), a cleaving plane (represented in black) passing through the points \(L_1(\Bx_1)\), \(L_2(\Bx_2)\), and \(L_3(\Bx_3)\), intersects \(\phi_0(\Bx)\) to create a region displayed in black in (d) where an obstacle can be placed. We roll, then, the cleaving plane about the line connecting \(L_1(\Bx_1)\) and \(L_3(\Bx_3)\), until the plane touches the point \(L_4(\Bx_4)\), see panels (e) and (f). At this point, we let the plane roll about the line connecting \(L_1(\Bx_1)\) and \(L_4(\Bx_4)\), until it reaches the point \(L_6(\Bx_6)\). Lastly, we let the cleaving plane roll about the line connecting \(L_1(\Bx_1)\) and \(L_6(\Bx_6)\) until it touches \(L_7(\Bx_7)\), as showed in (g) and (h). The rolling of the plane is illustrated step by step in  Movie S3 in the supplementary material. Note that, by applying this procedure to the remaining points \(\Bx_i\), \(i=2,..,7\), we can find the maximal feasible areas that an obstacle can occupy.}
	\labfig{Airy_7points}
\end{figure}

\cn
\subsection{Construction of strut nets avoiding obstacles}\labsect{Subsec_Construction_strut_nets}
\subsubsection{A strut net avoiding triangular and square inclusions} 
Let us search for a strut net that supports  a single active vertical force (pointing downward) applied to the point with coordinates $(0,4)$ in a given Cartesian plane. We examine  two supporting segments delimited by the points  $(-4,0)$ and $(3,4)$  (segment 1), and the points $(3,0)$ and $(4,0)$ (segment 2), which are placed on the opposite extremities of the strut net. Each of such segments is composed of 100 fixed reaction points uniformly spaced.
The searched strut net needs to avoid four obstacles, which are formed by three triangles \cb ($T1$, $T2$ and $T3$) and one rectangle ($R$). \cn The coordinates of the vertices of such polygons are defined as follows

\begin{itemize}
    \item \cb $T1$ \cn $[(1, 2);(3, 2);(3, 3)]$;
   \cb $T2$ \cn  $[(0, 1.3);  (0.2, 3.5); (-0.2, 3.5)]$;
   \cb $T3$ \cn  $[(-2.5, 0.5); (-1, 0.5);  (- 1, 2)]$;
    \item $R$ \cn  $[(-2, 1.2); (-2, 3); (-3, 3); (-3, 1.2)]$.
\end{itemize}

\fig{figpierre} illustrates four strut nets obtained through the linear programming algorithm of Section \sect{4} using different objective functions. The result in \fig{figpierre}(a) was obtained by 
minimizing the total weight of the structure, i.e., the integral over the boundary of the domain of the normal derivative of the strut net function. In this case, the optimal strut net is composed of 7 struts and 2 elementary loops. 
The strut nets in \fig{figpierre}(b,c,d) 
were instead obtained 
by minimizing the heights of the cleaving planes 
of the obstacles $T1$, $T2$ and $T3$ at their centers of mass, respectively. \cn
 By minimizing one of such functions, the goal is to enlarge the room available to the corresponding obstacle. It does not
always achieve this goal because the minimization can tilt the cleaving plane, or move the surrounding facets of the Airy
stress function. \cn 
The strut net in \fig{figpierre}(b) is again composed of 7 struts and 2 elementary loops, as that in (a). The strut net in (c) is instead formed by five struts and one loop, while that in (d) shows ten struts and two loops. \cn Observe that in most directions more
space is available to obstacle 3 in (c) rather than (d). Solutions like those in (c), where there are no struts separating two obstables,
can be enforced, if possible, by treating the two obstacles as a single obstacle.
Also, \cn it is worth noting that 
one could remove the last fork on the left of the strut net in \fig{figpierre}(d) and replace it by a straight line that touches the support midway. This is another way to \cn  enlarge \cn the room around T3 (not shown in the figure). \cn

\begin{figure}[h!]
	\centering
	\includegraphics[width=1.05\textwidth]{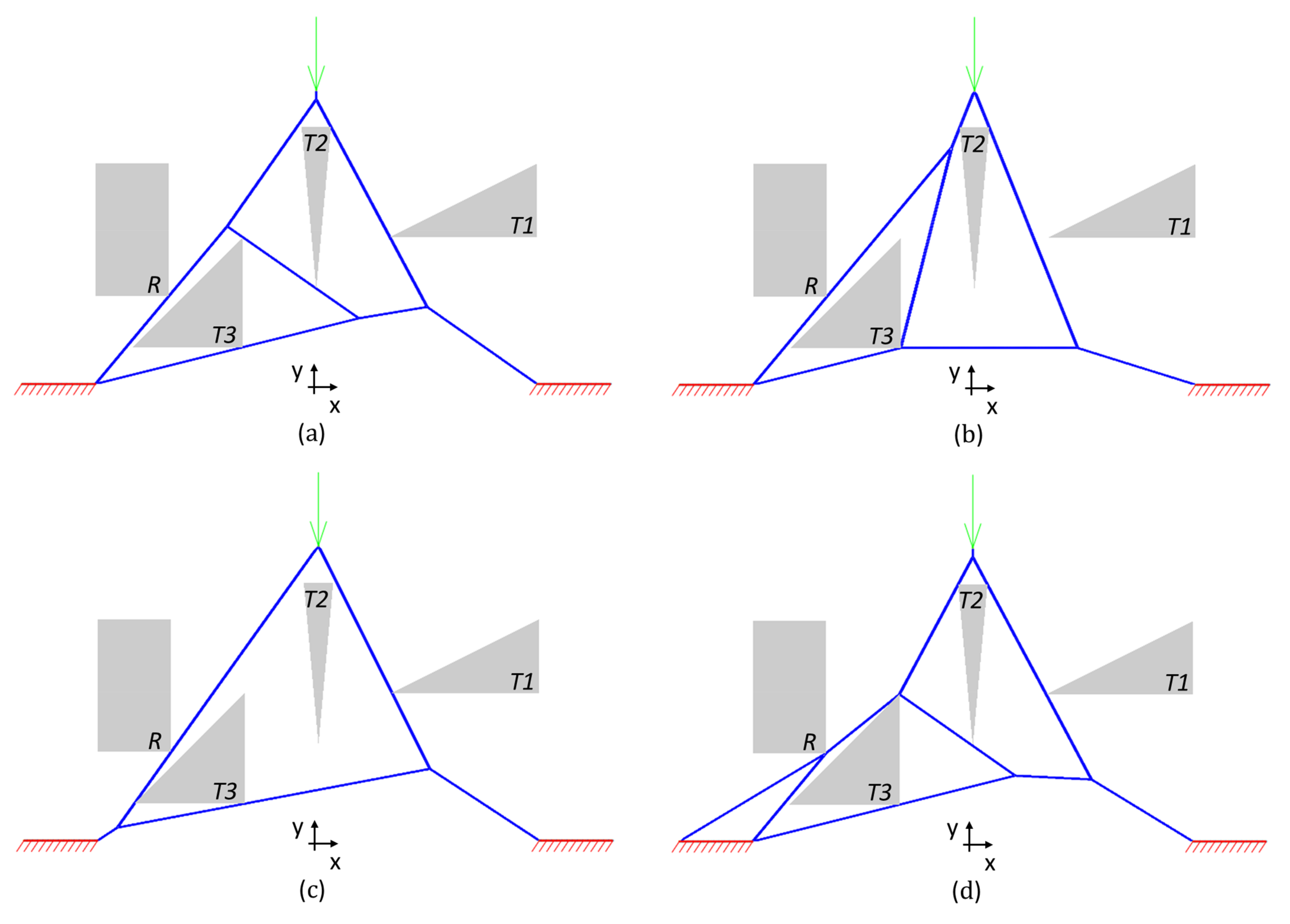}
	\caption{Strut nets resulting from the linear programming problem: a single active vertical force (marked in green) is supported by the shown strut net (blue in color) with the help of 200 reacting points ( placed along the segments marked in red).  The net must avoid 4 regions (marked in gray). The strut net in (a) minimizes the total weight of the structure, while those in (b,c,d) minimize the heights of the cleaving planes  of the obstacles $T1$, $T2$ and $T3$ at their centers of mass, respectively. \cn}
	\labfig{figpierre}
\end{figure}

\subsubsection{Arched strut nets avoiding elliptical inclusions}\labsect{arch2} 

Let us now pass to examine strut net models of masonry arches loaded by
vertical forces acting downwards on top of the arch. 
We consider \cg two \cn different sets of 23 vertical active forces, equally
spaced along the \cg horizontal \cn line segment connecting the point $(-2.75,4.20)$ with the
point $(2.75,4.20)$ of a given Cartesian plane (`loading segment' \crevv running along the extrados of the arch) \cn. We also examine two supporting segments, which coincide with those analyzed in the previous example.
The first obstacle to be avoided by the strut net is a half-disk of
radius $R=2.9$, centered at the origin of the Cartesian plane. It
represents the central region that the arch needs to cross (`arch obstacle').
The second and third obstacles are two elliptical regions, which are
respectively centered at the point $(1.2, 3.5)$ and the point $(-1.2,
3.5)$.
These obstacles exhibit vertical semi-major axis equal to $0.40$ and
horizontal semi-minor axis equal to $0.25$.
They are intended to reproduce holes that need to be drilled within the masonry arch  (`hole obstacles'). 
The masonry structure that contains the strut nets as `thrust lines' (or internal resisting structures) \cite{heyman1995} are colored light \crevv brown \cn in the figures that follow. It is composed of two triangular abutments with 1.10 width at the base and 4.20 height, and a filling region above the arch. All the examined strut nets lie within the masonry, which implies that the structure is stable under the examined loading conditions, according to the master safe theorem of masonry arches \cite{heyman1995}.
\crevv All the strut nets of the examples that follow have been obtained through the linear programming algorithm of Section \sect{4}, by
minimizing the total weight of the structure. \cn \cg In the simulations the tension in the outermost struts is very small and these struts disappear if the objective
function to be minimized is suitably perturbed.\cn

\begin{figure}[h!]
	\centering
	\includegraphics[width=0.65\textwidth]{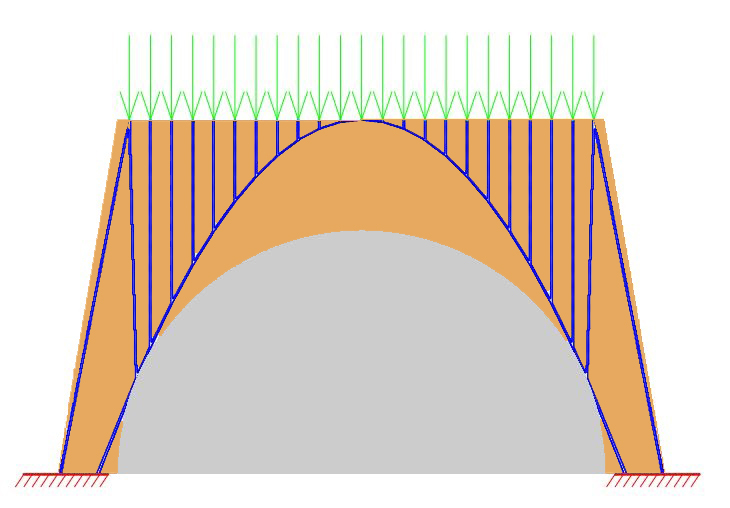}
	\caption{An arch that avoids a half disk (represented in gray), 
supporting 23 vertical forces of equal magnitude acting downwards on the top (in green) with the help of two 
supporting segments (in red). A feasible strut net (in dark blue) is obtained 
by using our linear programming formulation where each supporting 
segment has been replaced by 100 fixed reactive points, and the half disk
is approximated by a 101 sided polygon.
\crevv The masonry arch containing the examined strut net has been colored light brown. \cn
}
\label{arch0}
      \end{figure}

\begin{figure}[h!]
	\centering
	\includegraphics[width=0.65\textwidth]{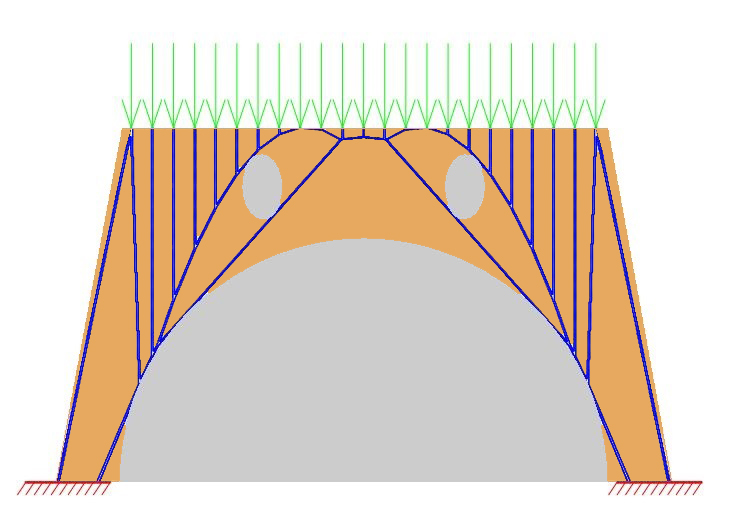}
	\caption{An arched strut net avoiding three
          obstacles. Like the previous example, it supports 23
          vertical downward forces of equal magnitude.  Additionally, the 
structure avoids the two elliptical regions, each approximated by a 20 sided
polygon.
}
\label{arch1}
      \end{figure}

We begin by considering the case when the 23  active vertical forces
are of equal (unit) magnitude and there is only the arch obstacle. The
strut net obtained for this set of
forces is illustrated in Figure \ref{arch0}. It consists of an open strut net formed by 49 struts, which is symmetric with respect to the $y$ axis. 
As a second example, we again consider 23 vertical active forces of
equal magnitude, as in the previous case, but this time we search for
a strut net that avoids both the arch obstacle and the two ellipse
obstacles defined above. Figure \ref{arch1} shows the solution obtained for the current example, which is composed of 53 struts, again symmetrically distributed with respect to the $y$ axis. The elliptical obstacles are embraced by closed loops of the strut net.


\newpage

\begin{figure} [h!]
	\centering
	\includegraphics[width=0.7\textwidth]{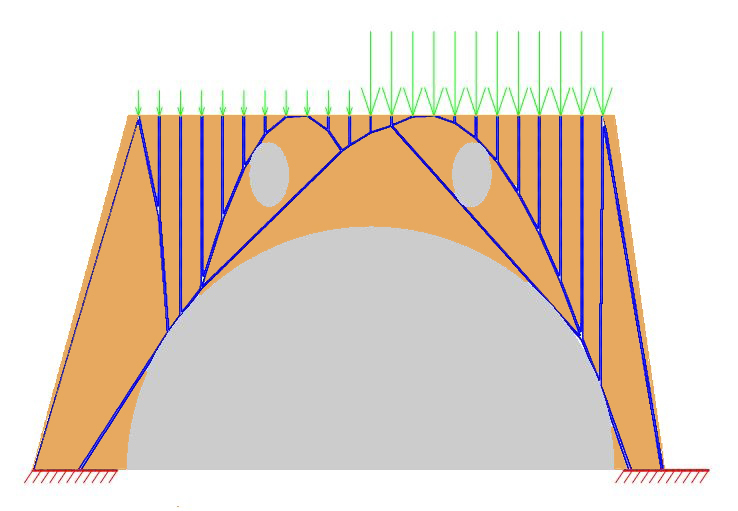}
	\caption{An arched strut net (in dark blue) avoiding three
          obstacles (in grey), which supports 11 vertical downward
          forces with magnitude 0.3, and 12 vertical downward forces with magnitude 1.0.
          }
\label{arch2}
      \end{figure}

We now pass to analyze a system of active forces composed of 11 vertical forces with magnitude $0.3$ uniformly spaced on the half of the loading segments placed along the negative $x$-axis; and 12 vertical forces with unit magnitude equally spaced along the complementary half of the loading segment.
The strut net returned by the linear programming algorithm for the
present case is shown in Figure \ref{arch2}. It is composed of 52 struts that are not symmetrically distributed with respect to the $y$ axis, and exhibits three closed loops. 

\section{Conclusions} \labsect{6}
     \setcounter{equation}{0}

Given a  two-dimensional \cn set of forces at the vertices of a convex polygon, some of
which are prescribed forces and the remainder reactive forces, our main
goal was to determine if a strut net exists avoiding a given set of
obstacles, and if so, to construct one such strut net. By approximating
each obstacle by a convex polygon, possibly with many sides,
we successfully devised an algorithm for doing this. It reduced to a
linear programming procedure and was based on finding a suitable
concave polyhedral Airy stress function associated with the strut
net when all struts are under compression. Additionally, under the
assumption that there is just one obstacle, we devised an algorithm
for generating regions that provide the maximum amount of space \cb available \cn to the
obstacle, \crev as \cn\cb shown \cb\crev in \fig{Airy_7points} and Movie S3. \cn By this we mean that the region available to the
obstacle cannot be enlarged in the sense of Definitions \ref{def1} and
\ref{def2}. We also obtained partial results (Theorem 1) in the case where there
are $q$ forces strictly inside the convex hull of all points
where forces are applied. We established that a strut net supporting these forces,
while avoiding $p$ obstacles, $p_0$ of which intersect this convex hull, can be replaced by a strut net, supporting
the same forces and avoiding the same obstacles, that has at most $q+p-p_0$ elementary loops. \crev In Figure \ref{arch_reduction} and Movie S1, $q=0$, and $p=p_0=1$ and the \cn\cb number of loops can be reduced until one obtains a final structure that is \cn\crev an open strut net. In \fig{loop_reduction} and Movie S2, instead, \cn\cb $q=p=1$ \cn\crev, and $p_0=0$, and the final net presents two elementary loops. \cn

Our work can be used to determine if a given masonry structure \crevv (modeled as \cn\cg an incompressible, \cn\crevv no-tension body \cite{heyman1995}) \cn can
support,  under compression, one or more families of given forces (and convex combinations of them).
For example,  the masonry arch colored light \crevv brown \cn in Figures
\ref{arch1} and \ref{arch2} can support at least 2 families of forces
(and convex combinations of them), namely those in the figures, while avoiding the two elliptical
obstacles. Alternatively, we may use the strut nets that support the
desired families of forces to efficiently design the masonry structures
themselves, eliminating unnecessary regions. Constraints on
the boundary of the masonry structure can be imposed by introducing
appropriate obstacles such as, for instance, the semicircular grey obstacle
beneath the arch in Figures \ref{arch0}, \ref{arch1} and \ref{arch2}.\cg Our analysis also
applies to structures built from unreinforced concrete since it supports compression but not tension. \cn

We address the following generalization of the present study \cn
 in \cn  future work:  given a force set
$\Bt=(\Bt_1,\Bt_2,\ldots,\Bt_n)$ and a vector set
$\Bf=(\Bf_1,\Bf_2,\ldots,\Bf_n)$,  \cn  determine \cn  the extreme values of
$\Gl$ such that a strut net avoiding the obstacles supports both $\Bt$
and $\Bt+\Gl \Bf$. Here $\Bf$ could represent an additional forcing,
say, e.g., due to an earthquake \cite{como2017,orduna17}.
\cn

There are many other avenues that have not been explored  in the present paper. \cn An extension
of the results to three dimensions would obviously be important, but
perhaps very difficult. The stress does have a matrix valued
function, the Beltrami stress function, that is analogous to the
two-dimensional Airy stress function. However, it is difficult to interpret the
constraints on the Beltrami stress function
imposed by a negative semidefinite stress field. Another
extension is to find a way of generating strut nets supporting $q$ forces strictly inside the convex hull of all points
where forces are applied, while avoiding $p$ obstacles, without
assuming one starts with a strut net achieving these goals. \cb Further goals are outlined in \cite{Bouchitte:2020:GSC}.
One, \cn with or without obstacles, is to allow for some elasticity in the struts, and to determine
the possible sets of displacements $\Bu_1,\Bu_2,\ldots,\Bu_n$  at the
vertices $\Bx_1,\Bx_2,\ldots,\Bx_n$ when a force set $\Bf_1,\Bf_2,\ldots,\Bf_n$
is applied. (Assuming the force set is such that at least one strut
net supporting the forces, yet avoiding any obstacles,  exists.) \cb Another \cn
is to allow for finite deformations, with or without obstacles.
That is, given moving points  $\Bx_1(t),\Bx_2(t),\ldots,\Bx_n(t)$ and
balanced forces  $\Bf_1(t),\Bf_2(t),\ldots,\Bf_n(t)$ that are
dependent on time $t$, when does there exist a single strut net under
tension that supports these forces and avoids any obstacles for some
given interval of $t$?  Here, single strut net means a strut net wherein the angles between struts change with $t$
but the topology and strut lengths do not; all struts also remain under compression and do not collide. \cb An even more challenging problem results \cn if some of the struts collide as $t$
changes.

\vskip6pt

\enlargethispage{20pt}





\section*{Acknowledgments}

OM and GWM are grateful to the National Science Foundation for support through Research Grants DMS-2008105 and DMS-2107926.
AA is grateful to the Italian Ministry of University and Research  for support through the PRIN 2017 grant 2017J4EAYB.

Fernando Fraternali is thanked for suggesting the incorporation of reactive forces and for suggesting many pertinent
  references. 

\section*{Data availability}

\crev Three movies are provided as supplementary material. Movie S1 illustrates the step-by-step loop reduction procedure for the example depicted in Figure \ref{arch_reduction}, whereas Movie S2 treats the loop reduction of the strut net represented in \fig{loop_reduction}. Finally, Movie S3 illustrates the rolling of the cleaving plane to create maximal regions where the obstacle can be placed, as illustrated in \fig{Airy_7points}.
\cn

\section*{Supplementary materials}
Movie\_S1.mp4

\medskip

\noindent \textbf{Caption for Movie S1}

\medskip

\noindent An arched strut net featuring $\ell=13$ closed loops, which
avoids an semi-circular arch obstacle ($p=1$),
is reduced to an open strut net ($\ell=0$) by employing the loop reduction procedure presented in \cite{Bouchitte:2019:FCW}.
The active and reactive external forces are at the vertices of
a convex polygon ($q=0$) \cn  and the obstacle intersects this polygon ($p_0=1$).
It is worth observing that the final arch satisfies $\ell\leq
q+p-p_0$, in agreement with Theorem 1 of Section \sect{2}.  
\noindent 

Movie\_S1.mp4

\medskip
\crev
\noindent \textbf{Caption for Movie S2}

\medskip
We consider an example in which one of the 5 applied forces lies inside the convex hull formed by the points where the forces are applied. A \cn\cb strut net \cb\crev  that supports such forces and avoids the obstacle represented by the gray circle is the one connecting the points pairwise. The movie shows how each elementary loop forming such a net can be replaced by an open \cn\cb strut net \cb\crev . In agreement with  Theorem 1 of  Section \sect{2}, the reduction can be done \cn\cb until \cb\crev at most 2 remaining elementary loops \cn\cb remain, \cb\crev as \cn\cb shown \cb\crev in the last frame of the movie. 

\noindent 

Movie\_S2.mp4

\medskip

\noindent \textbf{Caption for Movie S3}

\medskip
We consider 7 forces applied at the vertices of a convex polygon and we use the algorithm proposed at the end of Section \sect{3} to determine the maximal regions in which an obstacle can be placed so that a \cn\cb strut net \cb\crev  avoiding it exists. The movie starts with a picture of the open strut net function on the left-hand side, and the associated open \cn\cb strut net \cb\crev on the right-hand side. Secondly, the cleaving plane passing through the points \(L_1(\Bx_1)\), \(L_2(\Bx_2)\), and \(L_3(\Bx_3)\) is depicted on the left-hand panel, and the corresponding region \(\GG\) is depicted in black in the right-hand panel. Pivoting around point \(L_1(\Bx_1)\), the cleaving plane rolls about the line connecting \(L_1(\Bx_1)\) and the last point touched by the plane (initially \(L_3(\Bx_3)\)), spanning over the maximal elements \(\GG\). 
\noindent 

Movie\_S3.mp4
\cn



\end{document}